\documentclass[11pt]{article}
\usepackage{xcolor}
\usepackage{amsmath,amssymb,fullpage,color,url}
\usepackage[tmargin=1.0in,bmargin=1.0in,rmargin=1.0in,lmargin=1.0in]{geometry}
\usepackage[breaklinks=true]{hyperref}
\usepackage{graphicx}
\usepackage{verbatim}
\newtheorem{fact}[equation]{Fact}
\newtheorem{lemma}[subsubsection]{Lemma}
\newtheorem{theorem}[subsubsection]{Theorem}

\newtheorem{remark}[equation]{Remark}

{\medskip}

\newenvironment{proof}{\vspace{-0.05in}\noindent{\bf Proof.}}%
{\hspace*{\fill}$\Box$\par}
{\hspace*{\fill}$\Box$\par\vspace{4mm}}
{\hspace*{\fill}$\Box$\par}

\makeatletter
\newcommand{\mypm}{\mathbin{\mathpalette\@mypm\relax}}
\newcommand{\@mypm}[2]{\ooalign{%
		\raisebox{.1\height}{$#1+$}\cr
		\smash{\raisebox{-.6\height}{$#1-$}}\cr}}
\makeatother
\usepackage{blindtext}
\usepackage{textcomp}
\usepackage{commath}
\usepackage{mathtools}
\usepackage{amsmath}
\usepackage{mathabx}
\usepackage[colorinlistoftodos]{todonotes}
\usepackage{hyperref}

\title{Coupled K\"ahler-Einstein and Hermitian-Yang-Mills equations}
\author{Kartick Ghosh}
\date{}
\begin{document} 	 	
	\maketitle
	\begin{abstract} We introduce a new system of equations coupling K\"ahler-Einstein and Hermitian-Yang-Mills equations.  We provide a moment map interpretation of these equations. We identify a Futaki type invariant as an obstruction to the existence of solutions to these equations. We also prove a Matsushima-Lichnerowicz type theorem. We prove a deformation result that produces nontrivial solutions of these equations under some conditions. We produce examples on some projective bundles using Calabi ansatz. 
	\end{abstract}
\section{Introduction}
One of the main goals of differential geometry is to find ``nice" metrics on manifolds. In K\"ahler geometry, K\"ahler-Einstein metric is one such example. K\"ahler-Einstein metrics exist when the first Chern class of the manifold is negative (\cite{Aubin 1},\cite{Yau}) and when the first Chern class of the manifold is zero \cite{Yau}. When the first Chern class of the manifold is positive(the Fano case), there are obstructions. It is now known(\cite{CDS 1},\cite{CDS 2},\cite{CDS 3},\cite{Tian},\cite{Ved-Gabor}, \cite{Berman},\cite{Berman et al}, \cite{Chen et al}, \cite{Zhang} and see references therein) that such metrics exist if and only if an algebro-geometric obstruction called K-stability is met.\par
On the other hand, in the vector bundle setting, the goal is to find canonical connections. One such example is the Hermitian-Yang-Mills/ Hermitian-Einstein connection on holomorphic vector bundles over compact K\"ahler manifolds. The main result in this setting is the Kobayashi-Hitchin-Donaldson-Uhlenbeck-Yau correspondence that states that a holomorphic vector bundle over a compact K\"ahler manifold admits a Hermitian-Yang-Mills connection if and only if it is Mumford-polystable(\cite{Donaldson 1}, \cite{Donaldson 2}, \cite{Uhlenbeck and Yau}). The correspondence also holds true in the non K\"ahler case (\cite{Buchdahl}, \cite{Li And yau}).\par
The study of systems of equations coupling metrics on manifolds with connections on vector bundles is an active area in complex differential geometry. Such equations have been studied for a long time in the context of string theory in physics ( \cite{Li and Yau} for instance). Another important motivation to study coupled equations comes from their relationship with algebraic geometry, in particular, with the moduli problem for tuples consisting of polarised manifolds and holomorphic bundles over them.\par
In this paper, we couple these two fundamental equations and study them.
Let $(X, \omega)$ be a compact K\"ahler $n$-manifold with $-\sqrt{-1}\omega$ as the curvature form of the Chern connection on the Hermitian line bundle $(K_{X}^{-1},h)$. We can view the metric $h$ as a volume form $\Omega_{h}$.
Let $\pi: E \rightarrow X$ be a holomorphic vector bundle and  $H$ be a Hermitian metric on the bundle. 
In this paper, we consider the following coupled equations:
\begin{equation}
	\label{Coupled equations}
	\begin{split}
		&\sqrt{-1}{\bigwedge}_{\omega}F_{H}=\lambda Id\\
		&\frac{\alpha_{0}}{2}\left(\frac{\Omega_{h}}{\int_{X}\Omega_{h}}-\frac{\omega^{n}}{\int_{X}\omega^{n}}\right)-\alpha_{1}\frac{2}{(n-2)!}tr(F_{H}\wedge F_{H})\wedge\omega^{n-2}=\tilde{C}vol_{\omega},
	\end{split}
\end{equation}
where $F_{H}$ is the curvature of the Chern connection $A_{H}$ of $H$ on $E$; $\alpha_{0},\alpha_{1}\in \mathbb{R}$ are constants and $\tilde{C}$, $\lambda$ are topological constants. The topological constants are given as  $\tilde{C}=-\alpha_{1}\frac{\frac{2}{(n-2)!}\int_{X}tr(F_{H}\wedge F_{H})\wedge \omega^{n-2}}{\int_{X}^{}vol_{\omega}}$ and $\lambda=\frac{2\pi n}{n!vol_{\omega}(X)}\mu(E)$ where $\mu(E)=\frac{deg(E)}{rank(E)}$ is the slope of the vector bundle. Here the unknowns are a metric $h$ on $K_{X}^{-1}$ and a connection $A_{H}$ on $E$. \par
We provide a moment map interpretation of these equations following the method of \'Alvarez-C\'onsul--Garc\'ia-Fern\'andez--Garc\'ia-Prada \cite{Garcia et al} and using the moment map interpretation of K\"ahler-Einstein equation on Fano manifolds \cite{Donaldson 3}. This is the content of section (\ref{Moment Map Interpretation section}).\par
 Futaki \cite{Futaki} introduced an invariant which is an obstruction to the existence of K\"ahler-Einstein metric. The moment map interpretation of the K\"ahler-Einstein equation in the Fano case \cite{Donaldson 3} gives us an invariant. We prove that this invariant is the same as the classical Futaki invariant.
\begin{lemma}
	Suppose $(X,\omega)$ be a Fano manifold and $f$ is the Ricci potential i.e., $Ric(\omega)-\omega=\sqrt{-1}\partial\bar{\partial} f$ normalised so that $\int_{X}^{}e^{f}\omega^{n}=\int_{X}^{}\omega^{n}=C$. Then 
	\begin{equation}
		\int_{X} v(f)\omega^{n}=	\int_{X}\theta_{v}S(\omega)\omega^{n}=C\left(\int_{X}\theta_{v}\left(\frac{\Omega_{h}}{\int_{X}\Omega_{h}}-\frac{\omega^{n}}{C}\right)\right),
	\end{equation} 
	where $\theta_{v}$ is the Hamiltonian function corresponding to a holomorphic vector field $v$(i.e., $i_{v}\omega=\bar{\partial}\theta_{v}$) and $h$ is a Hermitian metric on $K_{X}^{-1}$ whose curvature is $\omega$. We normalize $\theta_{v}$ by requiring $\int_{X}\theta_{v}\omega^{n}=0$.
\end{lemma} 
Following the work of \'Alvarez-C\'onsul--Garc\'ia-Fern\'andez--Garc\'ia-Prada \cite{Garcia et al} for K\"ahler-Yang-Mills equations, we identify a Futaki type invariant which is an obstruction to the existence of solutions of the coupled equations (\ref{Coupled equations}). This is the content of section (\ref{Futaki invariant section}). \\
In section (\ref{matsushima obstruction section}), we introduce another obstruction which is analogous to the theorem of Matsushima-Lichnerowicz(\cite{Matsushima}, \cite{Lich 1}).\par
There are some trivial examples of solutions of our equations. On the Riemann sphere, we can produce examples of solutions considering polystable bundles over it since $F_{H}\wedge F_{H}$ vanishes for dimensional reason (\cite{Donaldson 4}, \cite{Narasimhan and Sheshadri}). For higher dimensional manifold, determining whether the coupled equations (\ref{Coupled equations}) admits solutions is a difficult problem because in this case the equations are coupled equations. If we consider a Fano manifold that admits a K\"ahler-Einstein metric and a holomorphic vector bundle over it which is polystable, then this serves as an example when $\alpha_{1}=0$.\par
So it is natural to wonder whether there are non-trivial examples. In section (\ref{Deformation section}) , we found non-trivial examples using deformation. In particular,  we prove the following deformation result.
\begin{theorem}
	Suppose $\omega$ is a K\"ahler-Einstein metric on a Fano manifold $X$ and $A$ is a Hermitian-Yang-Mills connection on a Hermitian holomorphic vector bundle $E$. Further, assume that $X$ has no non-zero holomorphic vector field. Then there exists $\epsilon>0$ such that for $\tilde{\alpha}\in\mathbb{R}$ with $-\epsilon<\tilde{\alpha}<\epsilon$, there is a solution $(\omega_{\phi},A_{\xi})$ to the coupled equation (\ref{Coupled equations}) with coupling constant $(1,\tilde{\alpha})$. 
\end{theorem} 
In our quest to find more nontrivial examples, we follow the method of Keller-T{\o}nnesen-Friedman  \cite{Keller} to produce non-trivial examples on some  projective bundles using Calabi ansatz (\cite{Calabi}, \cite{Hwang and Singer}, \cite{Tonnesen}, \cite{Gabor}). We state two theorems. The theorems look similar but the line bundles $L_{k}$ involved are different(for details see section \ref{Calabi ansatz section}).
\begin{theorem}
	Let $\Sigma_{i}=\mathbb{C}P^{1}$ and consider the line bundle $L=\otimes_{i=1}^{k}\pi_{i}^{*}(L_{i})$ over $\prod_{i=1}^{k}\Sigma_{i}$, where $L_{i}$ is a holomorphic line bundle of degree $-1$ over $\Sigma_{i}$. Set $X_{k}:=\mathbb{P}(L \oplus \mathcal{O})\rightarrow \prod_{i=1}^{k}\Sigma_{i}$, where $\mathcal{O}$ is the trivial line bundle over $\prod_{i=1}^{k}\Sigma_{i}$. Then there exist a metric $\omega\in 2\pi c_{1}(X_{k})$ and a connection $A_{H}$ on a line bundle $L_{k}$ such that they satisfy the coupled equations (\ref{Coupled equations}) if the following conditions are met:
	\begin{enumerate}
		\item 	Case 1: $k=1$\\
		$(1) \alpha_{0},\alpha_{1}>0$\\
		$(2) \alpha_{0}=8b_{1}^{2}\alpha_{1}(13\log3-12)$
		\item	Case 2: $k=2,3,4$\\
		$(1) \alpha_{0},\alpha_{1}>0$\\
		$(2) \frac{k+1}{2(2\pi)^{k+1}(3^{k+1}-1)}\alpha_{0}=\alpha_{1}\frac{2b_{k}^{2}(k+1)k}{(2\pi)^{2}}\left(\frac{(k+1)(3^{k-1}-1)}{(k-1)(3^{k+1}-1)}-R_{k}\right)$,
	\end{enumerate}
	where $R_{k}=\frac{\int_{0}^{2}(1-t)(1+t)^{k-2}dt}{\int_{0}^{2}(1-t)(1+t)^{k}dt}$ and $b_{k}=\frac{2m_{2}-3m_{1}}{2+3k\log3}$ for some integers $(m_{1},m_{2})\neq(0,0)$. Moreover, for $k\geq 5$ our method does not produce any solutions.
\end{theorem}
 Now we state another theorem where the manifold is the same but the line bundle is different(for details see \ref{Calabi ansatz section}).
\begin{theorem}
	Let $\Sigma_{i}=\mathbb{C}P^{1}$ and consider the line bundle $L=\otimes_{i=1}^{k}\pi_{i}^{*}(L_{i})$ over $\prod_{i=1}^{k}\Sigma_{i}$, where $L_{i}$ is a holomorphic line bundle of degree $-1$ over $\Sigma_{i}$. Set $X_{k}:=\mathbb{P}(L \oplus \mathcal{O})\rightarrow \prod_{i=1}^{k}\Sigma_{i}$, where $\mathcal{O}$ is the trivial line bundle over $\prod_{i=1}^{k}\Sigma_{i}$. Then there exist a metric $\omega\in 2\pi c_{1}(X_{k})$ and a connection $A_{H}$ on a line bundle $L_{k}$ such that they satisfy the coupled equations (\ref{Coupled equations}) if the following conditions are met:
	\begin{enumerate}
		\item 	Case 1: $k=2$\\
	$(1) \alpha_{0},\alpha_{1}>0$\\
	$(2) \alpha_{0}=32\pi\alpha_{1}$
\item	Case 2: $k=3$\\
	$(1) \alpha_{0},\alpha_{1}>0$\\
	$(2) \frac{\alpha_{0}}{32(2\pi)^{2}}=\frac{38}{21}\alpha_{1}$
\item	Case 3: $k=4$\\
	$(1) \alpha_{0},\alpha_{1}>0$\\
	$(2) \frac{3\alpha_{0}}{4!(2\pi)^{3}}=\frac{2\times 235}{3\times 23}\alpha_{1}$.
\end{enumerate}
Moreover, for $k\geq 5$ our method does not produce any solutions.
\end{theorem}
We plan on exploring these equations further in future work.\\
\emph{Acknowledgments:} I thank my advisor, Vamsi Pritham Pingali, for suggesting this problem to me and for his constant encouragement. He helped me to correct several mistakes and make the paper more readable. The author is supported by a scholarship from the Indian Institute of Science. 
\section{Moment Map Interpretation}
\label{Moment Map Interpretation section}
In this section, we provide the moment map interpretation of the coupled equations. To achieve that we follow the method of \'Alvarez-C\'onsul--Garc\'ia-Fern\'andez--Garc\'ia-Prada \cite{Garcia et al}.
\subsection{Moment map interpretation of Hermitian-Yang-Mills equation}
\label{HYM moment map section}
Let $(X,\omega)$ be a compact symplectic manifold of dimension $2n$ and $E$ is a Hermitian vector bundle over it. Let $\mathcal{A}$ be the space of unitary connections on $E$.  The tangent bundle $T\mathcal{A}$ can be identified as $\Omega^{1}_{X}(T_{sh})$, where $T_{sh}\subset End(E)$ is a bundle of skew Hermitian endomorphisms. We endow $\mathcal{A}$ with the symplectic form $\omega_{\mathcal{A}}$ defined by 
\begin{equation}
\label{HYM symplectic form}
\omega_{\mathcal{A}}(a,b)=-\int_{X}tr(a\wedge b)\wedge \omega^{n-1}, 
\end{equation}
where $a,b \in T\mathcal{A}\cong \Omega^{1}_{X}(T_{sh})$. The gauge group of $E$ denoted by $\mathcal{G}$ is the group of automorphisms of $E$ covering the identity on $X$. This group $\mathcal{G}$ acts on $\mathcal{A}$. Atiyah and Bott \cite{Atiyah-Bott} when $X$ is a Riemann Surface and  Donaldson (\cite{Donaldson 1} , \cite{Donaldson 2}) in higher dimension proved that the action is Hamiltonian with equivariant moment map $\mu_{\mathcal{G}}:\mathcal{A} \mapsto (Lie\mathcal{G})^{*}$ given by
 \begin{equation}
 	\langle\mu_{\mathcal{G}}(A),\xi\rangle=\sqrt{-1}\int_{X}tr\left(\xi\wedge(\sqrt{-1}{\bigwedge}_{\omega}F_{A}-\lambda Id)\right)\frac{\omega^{n}}{n!},
 \end{equation}
where $\xi\in Lie\mathcal{G}$ and $\lambda$ is a topological constant.\par
Now suppose $(X,\omega,J)$ is a compact K\"ahler manifold.  Then there is a subspace  $\mathcal{A}^{1,1}$ of $\mathcal{A}$ consisting of connections with $F_{A}\in \Omega^{1,1}(End(E))$. This subspace is in one to one correspondence with holomorphic structures on $E$.  The complex structure $J$ induces a natural complex structure on $\mathcal{A}^{1,1}$ by 
\begin{equation}
	\label{almost complex structure HE}
	J_{\mathcal{A}}:=a \mapsto -a(J.).
\end{equation}
This complex structure is compatible with $\omega_{\mathcal{A}}$. This makes $(\mathcal{A}^{1,1},\omega_{\mathcal{A}})$ a K\"ahler manifold.
\subsection{Coupling method}
 We now briefly explain the general method of coupling(for more details see \cite{Garcia et al}). First, we describe the Hamiltonian action of an extended Lie group on a symplectic manifold. This action can be described under certain assumptions, in terms of a normal Lie subgroup and the quotient Lie subgroup. Then we apply this method where the symplectic manifold is the space of connections. 
\subsubsection{Extended Lie group and its action} 
\label{General coupling method section}
Suppose that there is an extension of Lie groups 
\begin{equation}
	\label{Lie group extension}
1\rightarrow\mathcal{G}\xrightarrow{\text{$i$}}\widetilde{\mathcal{G}}\xrightarrow{\text{$p$}}\mathcal{H}\rightarrow 1.
\end{equation}
This determines another extension 
\begin{equation}
\label{Lie algebra extension}
0\rightarrow Lie\mathcal{G}\xrightarrow{\text{$i$}}Lie\widetilde{\mathcal{G}}\xrightarrow{\text{$p$}}Lie\mathcal{H}\rightarrow 0.
\end{equation}
Let $\mathcal{A}$ be a manifold with an action of the extended Lie group $\widetilde{\mathcal{G}}$.\par
The short exact sequence (\ref{Lie algebra extension}) does not generally split as a sequence of Lie algebras, but it does always split as a sequence of vector spaces. Let $W \subset Hom(Lie \tilde{\mathcal{G}},Lie \mathcal{G})$ be the space of vector space splittings. Since $\mathcal{G}\subset \tilde{\mathcal{G}}$ is a normal subgroup,
 there is a well-defined $\tilde{\mathcal{G}}$-action
on $W$, given by $g.\theta\coloneqq Ad(g)\circ \theta \circ Ad(g^{-1})$ for $g\in \tilde{\mathcal{G}},\theta \in W$.
Let $\mathcal{W}\subset C^{\infty}(\mathcal{A},W)$ be the space of $\tilde{\mathcal{G}}$-equivariant smooth maps $\theta: \mathcal{A} \rightarrow W $. We consider the case where $\mathcal{W}$ is nonempty. \par
Let $\omega_{\mathcal{A}}$ be a symplectic form on $\mathcal{A}$ such that the $\widetilde{\mathcal{G}}$ action is symplectic. The aim is now to characterise the Hamiltonian $\widetilde{\mathcal{G}}$-action on $\mathcal{A}$ in terms of $\mathcal{G}$ and $\mathcal{H}$ under the assumption of existence of some $\theta\in \mathcal{W}$. Suppose that the $\tilde{\mathcal{G}}$-action on $\mathcal{A}$
is Hamiltonian with a $\tilde{\mathcal{G}}$-equivariant moment map $\mu_{\tilde{\mathcal{G}}}: \mathcal{A} \rightarrow (Lie\tilde{\mathcal{G}})^{*}$. We will break this
map into pieces corresponding to the Lie algebras $Lie \mathcal{G}$ and $Lie \mathcal{H}$ using $\theta \in \mathcal{W}$. First note that the
map
\begin{equation}
W \rightarrow Hom(Lie \mathcal{H}, Lie\tilde{\mathcal{G}}):\theta \rightarrow \theta^{\perp},
\end{equation}
where $\theta^{\perp}: Lie \mathcal{H} \rightarrow Lie \tilde{\mathcal{G}}$ is uniquely defined by the equation
\begin{equation}
	\label{5}
Id -i \circ \theta  = \theta^{\perp} \circ p,
\end{equation}
is $\tilde{\mathcal{G}}$-equivariant
with respect to the action in $Hom(Lie\mathcal{H},Lie\tilde{\mathcal{G}})$ given by $g.\theta^{\perp}=Ad(g)\circ \theta^{\perp}\circ Ad(p(g^{-1}))$ for all $g \in \tilde{\mathcal{G}}.$ Now we can break the moment map into 
\begin{equation}
	\label{2}
	\langle \mu_{\tilde{\mathcal{G}}},\xi\rangle =\langle \mu_{\tilde{\mathcal{G}}},i\theta\xi\rangle+\langle\mu_{\tilde{\mathcal{G}}},\theta^{\perp}p(\xi)\rangle,
\end{equation}
for all $\xi \in Lie\tilde{\mathcal{G}}$, where the summands on the right hand side define a pair of $\tilde{\mathcal{G}}$-equivariant maps $\mu_{\mathcal{G}}: \mathcal{A}\rightarrow (Lie\mathcal{G})^{*}$ given by
\begin{equation}
	\label{3}
	\langle\mu_{\mathcal{G}},\zeta\rangle \coloneqq\langle \mu_{\tilde{\mathcal{G}}},i\zeta\rangle
\end{equation}
for all $\zeta\in Lie\mathcal{G}$, and $\sigma_{\theta}:\mathcal{A}\rightarrow(Lie\mathcal{H})^{*}$
\begin{equation}
	\label{4}
		\langle \sigma_{\theta},\eta\rangle \coloneqq\langle \mu_{\tilde{\mathcal{G}}},\theta^{\perp}\eta\rangle
\end{equation}
for all $\eta \in Lie\mathcal{H} $.
Given a smooth map $\xi:\mathcal{A}\rightarrow Lie\tilde{\mathcal{G}}$, let us denote by $Y_{\xi}$ the
vector field on $\mathcal{A}$ given by
$(Y_{\xi})_{A}\coloneqq \frac{d}{dt}\bigg|_{t=0}exp(t\xi_{A}).A$, for all $A \in \mathcal{A}$.
In particular, any smooth $\theta:\mathcal{A}\rightarrow W$ defines a map
$Y_{\theta^{\perp}}: Lie \mathcal{H} \rightarrow \Omega^{0}(T\mathcal{A})$ defined by $\eta \rightarrow Y_{\theta^{\perp}\eta}$. The following theorem is the general method of coupling ( proposition $1.3$ in \cite{Garcia et al})
\begin{theorem}
	\label{theorem iff}
 The $\tilde{\mathcal{G}}$-action
	on $\mathcal{A}$ is
	hamiltonian if and only if so is the action of $\mathcal{G}\subset \tilde{\mathcal{G}}$ on $\mathcal{A}$, with a $\tilde{\mathcal{G}}$-equivariant moment
	map $\mu_{\mathcal{G}}: \mathcal{A}\rightarrow (Lie\mathcal{G})^{*}$, and there exists a smooth $\tilde{\mathcal{G}}$-equivariant
	map $\sigma_{\theta}:\mathcal{A}\rightarrow(Lie\mathcal{H})^{*}$
	satisfying
	\begin{equation}
		\label{non spliting condition}		
		i_{Y_{\theta^{\perp}\eta}}\omega_{\mathcal{A}}=\langle\mu_{\mathcal{G}},d\theta(\eta) \rangle +d\langle \sigma_{\theta},\eta\rangle
	\end{equation}
for all $\eta \in Lie\mathcal{H}$. By definition, $d\theta$ is a $\tilde{\mathcal{G}}$-invariant $Hom(Lie\mathcal{H},Lie\mathcal{G})$-valued $1$-form on $\mathcal{A}$.
	In this case, a $\tilde{\mathcal{G}}$-equivariant
	moment map $\mu_{\tilde{\mathcal{G}}}:\mathcal{A}\rightarrow(Lie\tilde{\mathcal{G}})^{*}$ is given by 
	\begin{equation}
		\label{7}
		\langle \mu_{\tilde{\mathcal{G}}},\xi\rangle=\langle\mu_{\mathcal{G}},\theta\xi\rangle +\langle\sigma_{\theta},p(\xi)\rangle
	\end{equation}
for all $\xi \in Lie\tilde{\mathcal{G}}$.
\end{theorem}
One can see that condition (\ref{non spliting condition}) generalizes the usual infinitesimal condition $i_{Y_{\eta}}\omega_{\mathcal{A}}=d\langle\mu_{\mathcal{H}},\eta\rangle$ for moments map $\mu_{\mathcal{H}}$
 for the induced $\mathcal{H}$-action on $\mathcal{A}$ when the Lie group extension (\ref{Lie group extension}) splits. 
 
 \subsubsection{Extended gauge group and its action on the space of connections}
 We apply the theory developed in section (\ref{General coupling method section}) to compute the moment map for the action of the extended gauge group of a bundle over a compact symplectic manifold, on the space of connections.\\
 Let $(X,\omega)$ be a compact symplectic manifold of dimension $2n$ and $(E,H)$ be a rank $r$ Hermitian vector bundle over it. Let $\mathcal{A}$ be the space of unitary connections on $(E,H)$. Let $\mathcal{H}$ be the group of Hamiltonian symplectomorphisms. Let $E_{H}$ be the principal $U(r)$-bundle of unitary frames of $(E,H)$. We define the extended gauge group $\widetilde{\mathcal{G}}$ to be the automorphisms of $E_{H}$ which covers elements of $\mathcal{H}$. Then the gauge group $\mathcal{G}$, already defined in section (\ref{HYM moment map section}) , is the normal subgroup $\mathcal{G}\subset\widetilde{\mathcal{G}}$ of automorphisms covering the identity.\par
There is a canonical short exact sequence of Lie groups
 \begin{equation}
 	\label{Special Lie group extension}
 	1\rightarrow\mathcal{G}\xrightarrow{\text{$i$}}\widetilde{\mathcal{G}}\xrightarrow{\text{$p$}}\mathcal{H}\rightarrow 1,
 \end{equation}
where $i$ is the inclusion map and the map $p$ assigns to each automorphism $g$ the Hamiltonian symplectomorphism it covers.  For the surjectivity of the map $p$ (see \cite{Garcia et al}).\\
 To every vector bundle $\pi:E\rightarrow X$, we can associate an exact sequence of vector bundles:
 	\[0\rightarrow VE_{H} \rightarrow TE_{H} \rightarrow \pi^{*}TX \rightarrow 0,\]
 	where $VE_{H}\subset TE_{H}$ is the vertical bundle of $E_{H}$. Elements of $\mathcal{A}$ can be viewed as splittings $A:TE_{H} \rightarrow VE_{H}$ of the above exact sequence. Using this viewpoint, we can define the action of $\widetilde{\mathcal{G}}$ on $\mathcal{A}$. The action $g.A$ is defined as $g\circ A \circ g^{-1}$(here the action is in the sense of infinitesimal action). Any such splitting $A$ induces a vector space splitting of the Atiyah short exact sequence 
 	\begin{equation}
 		0\rightarrow Lie(\mathcal{G})\xrightarrow{\text{$i$}}Lie(Aut(E_{H}))\xrightarrow{\text{$p$}}Lie(Diff(X))\rightarrow 0,
 	\end{equation}
 where $Diff(X)$ is the Lie algebra of vector fields on $X$ and $Lie(Aut(E_{H}))$ is the Lie algebra of vector fields on $E_{H}$. The splitting is given by maps 
 \begin{equation}
 	\label{Splitting map}
\theta_{A}:Lie(Aut(E_{H}))\longmapsto Lie\mathcal{G}, \ \ \theta_{A}^{\perp}: Lie(Diff(X))\longmapsto Lie(Aut(E_{H})) 
\end{equation}
 satisfying $i\circ \theta_{A}+ \theta_{A}^{\perp}\circ p=Id$, where $\theta_{A}$ is the vertical projection given by $A$ and $\theta_{A}^{\perp}$ is the horizontal lift of vector fields on $X$ to vector fields on $E_{H}$ given by $A$.\\
 The splitting (\ref{Splitting map}) restricts to a splitting of the short exact sequence 
 	\begin{equation}
 		0\rightarrow Lie(\mathcal{G})\xrightarrow{\text{$i$}}Lie(\widetilde{\mathcal{G}})\xrightarrow{\text{$p$}}Lie(\mathcal{H})\rightarrow 0
 	\end{equation}
 induced by (\ref{Special Lie group extension}) . Following the notation of section (\ref{General coupling method section}) , we see that 
 \begin{equation}
 	\begin{split}
 		&\theta: \mathcal{A}\mapsto W\\
 		&A\mapsto \theta_{A}
 	\end{split}
 \end{equation}
 is a $\widetilde{\mathcal{G}}$-equivariant smooth map. It is also clear that the $\widetilde{\mathcal{G}}$-action on $\mathcal{A}$ is symplectic,
 for the symplectic form (\ref{HYM symplectic form}) . We apply here the methods of section  (\ref{General coupling method section}) to get the following theorem (proposition $1.6$ in \cite{Garcia et al}).
\begin{theorem}
	\label{Moment map in coupling method}
	The $\tilde{\mathcal{G}}$-action on $\mathcal{A}$ is hamiltonian, with $\tilde{\mathcal{G}}$-equivariant moment map $\mu_{\tilde{\mathcal{G}}}:\mathcal{A}\rightarrow(Lie\tilde{\mathcal{G}})^{*}$ given by
	\begin{equation}
		\label{10}
		\langle\mu_{\tilde{\mathcal{G}}},\xi\rangle=\langle\mu_{\mathcal{G}},\theta\xi\rangle+\langle\sigma,p(\xi)\rangle \  for \  all \ \xi \in Lie\tilde{\mathcal{G}},
	\end{equation}
where $\mu_{\mathcal{G}}: \mathcal{A} \rightarrow (Lie\mathcal{G})^{*}$ and $\sigma : \mathcal{A}\rightarrow (Lie\mathcal{H})^{*}$ are given by
\begin{equation}
	\label{11}
	\begin{split}
		&\langle\mu_{\mathcal{G}},\theta\xi\rangle(A)=\sqrt{-1}\int_{X} tr(\theta_{A}\xi\wedge (\sqrt{-1}{\bigwedge}_{\omega}F_{A}-\lambda Id))\wedge \frac{\omega^{n}}{n!}\\
	&\langle\sigma,\eta_{\phi}\rangle(A)=\frac{1}{4}\int_{X} \phi \left(\frac{2}{(n-2)!}tr(F_{A}\wedge F_{A})\wedge\omega^{n-2}+4\lambda\frac{1}{(n-1)!} tr(\sqrt{-1}F_{A})\wedge\omega^{n-1}\right)\\
\end{split}
\end{equation}
for any $\eta_{\phi} \in Lie\mathcal{H}$.
\end{theorem}

\subsection{Moment map interpretation of KE equation}
In this section, we briefly explain the moment map interpretation of K\"ahler-Einstein equation in the Fano case. Donaldson \cite{Donaldson 3} first came up with the interpretation. We briefly explain Donaldson's moment map.\par
Suppose $U$ is a real vector space of dimension $2n$. For $\alpha,\beta\in \Lambda^{n}(U)\otimes\mathbb{C}$, we write
\[\langle\alpha,\beta\rangle=(\sqrt{-1})^{n^{2}}\alpha\wedge\bar{\beta}.\]
  This is a $\Lambda^{2n}\otimes \mathbb{C}$-valued indefinite Hermitian form on $\Lambda^{n}(U)\otimes\mathbb{C}$. Suppose that $U^{*}$ has a symplectic structure $\omega\in\Lambda^{2}(U)$. Let $M$ be the set of compatible almost complex structures on $U$. Then $M$ can be identified with an open subset of the Lagrangian Grassmannian of $U\otimes\mathbb{C}$. Thus $M$ can be identified with a subset of $\mathbb{P}({\Lambda}^{n}U\otimes\mathbb{C})$. Let $N\subset {\Lambda}^{n}U\otimes\mathbb{C}$ be the corresponding cone, with the origin deleted.\par
We now apply the above vector space constructions in the bundle setting where cotangent space will serve as the vector space. Let $(X,\omega)$ be a compact symplectic manifold of dimension $2n$ and $L$ be a Hermitian line bundle with metric connection $\nabla$ having $-\sqrt{-1}\omega$ as curvature. We denote $\Omega^{n}(L)$ as the $L$-valued $n$-forms. We define a Hermitian form on $\Omega^{n}(L)$ by
\begin{equation}
	\label{Hermitian form KE}
	\omega_{L}(\alpha,\beta)=\int_{X}\langle\alpha,\beta\rangle,
\end{equation}
where  $\langle \alpha,\beta\rangle=(\sqrt{-1})^{n^{2}}\alpha\wedge\bar{\beta}$ and we have used the Hermitian form on $L$ for the section part. We write $d_{\nabla}$ for the coupled exterior derivative on $L$-valued forms. We have an integration-by-parts formula 
\begin{equation}
	\label{Integration by parts formula}
\int_{X}\langle d_{\nabla}\alpha,\beta\rangle:=(-1)^{n}\int_{X}\langle\alpha,d_{\nabla}\beta\rangle,
\end{equation}
for $\alpha\in \Omega^{n-1}(L), \beta\in\Omega^{n}(L)$, where the definition of $\langle\ ,\ \rangle$ is extended in the obvious way. \\
 Let $\widetilde{\mathcal{H}}$ be the group of automorphisms of $L$, as a line bundle with connection, covering symplectomorphisms of $X$. The group $\widetilde{\mathcal{H}}$ acts on $L$ and on $X$ and hence on $\Omega^{n}(L)$. We denote this action by $h\circ \alpha$, where $h\in\widetilde{\mathcal{H}}$ and $\alpha\in \Omega^{n}(L)$. The Lie algebra of this group is $C^{\infty}(X)$. The infinitesimal action of the Lie algebra on $\Omega^{n}(L)$ is given by 
\begin{equation}
	\label{Infinitesimal action in KE case}
	\sqrt{-1}R_{\alpha}(H)=v_{H}(\dashv d_{\nabla}\alpha)-d_{\nabla}(v_{H}\dashv\alpha)-\sqrt{-1}H\alpha.
\end{equation} 
 Let $\underline{M}\rightarrow X$ be the bundle with fibre $M$. Thus sections of $\underline{M}$ are almost-complex structures on $X$ compatible with the symplectic form. Now we restrict ourselves to the case when the first Chern class of $L$ is equal to the first Chern class of any compatible almost-complex structure. Now consider $\underline{N}\rightarrow X$ as the bundle with fibre $N$, coupled to $L$ as above. Hence sections of $\underline{N}$ are $L$-valued $n$-forms on $X$ and the Chern conditions ensure that the global sections of $\underline{N}$ exist. $\underline{N}$ can be thought of as a principal $\mathbb{C}^{*}$ bundle over $\underline{M}$. Let $\mathcal{J}_{int}$ be the space of integrable compatible almost complex structures. The following lemma (lemma $1$ in \cite{Donaldson 3}) states how an integrable almost complex structure is related to a section of $\underline{N}$ and vice versa. 

\begin{lemma}
	\label{lemma 1}
	Any section $\alpha$ of $\underline{N}$ with $d_{\nabla}\alpha=0$ projects to an integrable almost complex structure. Conversely, if $J\in \mathcal{J}_{int}$ is an integrable complex structure it has a lift to a section $\alpha$ of $\underline{N}$ with  $d_{\nabla}\alpha=0$ and $\alpha$ is unique up to multiplication by a
	constant in $\mathbb{C}^{*}$.
\end{lemma} 

Let $\hat{\mathcal{J}}_{int}$ be the sections $\alpha$ of $\hat{\underline{M}}$ with $d_{\nabla}\alpha=0$. Then (\ref{lemma 1}) tells us that $\hat{\mathcal{J}}_{int}$ is a $\mathbb{C}^{*}$-bundle over $\mathcal{J}_{int}$. Now $\hat{\mathcal{J}}_{int}$ is a subspace of the vector space $\Omega^{n}(L)$ which has a Hermitian form given by (\ref{Hermitian form KE}) . One can check that $\omega_{L}(\alpha,\alpha)\geq 0$ for any $\alpha\in \hat{\mathcal{J}}_{int}$. We now describe the tangent space of $\hat{\mathcal{J}}_{int}$ at $\alpha$. The tangent space $T_{\alpha}$ consists of forms $\beta$ with $d_{\nabla}\beta=0$ and at each point $x\in X$ lie in the corresponding vertical tangent space of $\underline{N}$. The next lemma asserts where the Hermitian form is definite.
\begin{lemma}
	\label{lemma 2}
	The Hermitian form $\omega_{L}$ is negative-deﬁnite on the orthogonal
	complement of $\alpha$ in $T_{\alpha}$.
\end{lemma}
The previous lemma asserts that $\widetilde{\omega_{L}}\coloneqq-\omega_{L}$ induces a K\"ahler metric on $\mathcal{J}_{int}$. We denote the underlying complex structure by $J_{KE}$.\par
The action of the group $\widetilde{\mathcal{H}}$ is Hamiltonian and the moment map is given by 
\[\alpha\mapsto -\frac{1}{2}(H\rightarrow Re\ \omega_{L}(\alpha,R_{\alpha}H)),\]
where $H$ is a Lie algebra element. If we restrict our forms to the forms $\alpha$ with  $d_{\nabla}\alpha=0$, then we see that $\omega_{L}(\alpha,R_{\alpha}H)=\omega_{L}(\alpha,H\alpha)$(using \ref{Integration by parts formula}). Thus the moment map $\mu_{0}$, for the action on $\hat{\mathcal{J}}_{int}$ is given by \[\mu_{0}(\alpha)=-\frac{1}{2}\langle\alpha,\alpha\rangle,\]
using the pairing between $2n$-forms and functions. Since the action of $\widetilde{\mathcal{H}}$ commutes with the action of $S^{1}$ action, the moment map on the symplectic quotient $\mathcal{J}_{int}$ is given by \[\mu_{0}([\alpha])=-\frac{1}{2}\frac{\langle\alpha,\alpha\rangle}{\omega_{L}(\alpha,\alpha)}.\]
Adjusting the moment map by a constant, we can write 
\begin{equation}
	\label{moment map in KE not in final form}
	\mu([\alpha])=-\frac{1}{2}\left(\frac{\langle\alpha,\alpha\rangle}{\omega_{L}(\alpha,\alpha)}-\frac{\omega^{n}}{\int_{X}\omega^{n}}\right).
\end{equation}  
Suppose $X_{\alpha}$ is the complex manifold corresponding to $[\alpha]$. Using this $\alpha$, one can identify $L$ with $K_{X}^{-1}$ and therefore we get a metric $h$ on $K_{X}^{-1}$. This metric can be viewed as a volume form $\Omega_{h}$. In local coordinates $\Omega_{h}$ takes the form \[\Omega_{h}=V(\sqrt{-1})^{n}dz^{1}\wedge d\bar{z}^{1}\dots\wedge dz^{n}\wedge d\bar{z}^{n},\]
where $V=\lvert \frac{\partial}{\partial z^{1}}\wedge\dots\wedge\frac{\partial}{\partial z^{n}}\rvert_{h}^{2}$. In this setting, we have $\Omega_{h}=\langle\alpha,\alpha\rangle$. Now we can write the moment map (\ref{moment map in KE not in final form}) as 
\begin{equation}
	\label{moment map for KE equation}
	\mu([\alpha])=-\frac{1}{2}\left(\frac{\Omega_{h}}{\int_{X}\Omega_{h}}-\frac{\omega^{n}}{\int_{X}\omega^{n}}\right).
\end{equation}
\subsection{The coupled equation as a moment map condition}
\label{Coupled moment map section}
Again as before, let $(X,\omega)$ be a symplectic manifold of dimension $2n$ and $E$ be a Hermitian vector bundle over it. We consider the product manifold $\mathcal{A}\times \Omega^{n}(L)$ with the symplectic form 
\begin{equation}
	\label{Coupled symplectic form}
\omega_{\alpha}=\alpha_{0}\widetilde{\omega_{L}} +4\alpha_{1}\omega_{\mathcal{A}},
\end{equation}
 for pair of non-zero real constants $\alpha=(\alpha_{0},\alpha_{1})$. Using the complex structure $J_{\mathcal{A}}$(\ref{almost complex structure HE}) on $\mathcal{A}^{1,1}$ and $J_{KE}$ on $\mathcal{J}_{int}$, we can give a formally integrable almost complex structure $J$ on $\mathcal{A}^{1,1}\times \mathcal{J}_{int}$ which is compatible with $\omega_{\alpha}$ (\ref{Coupled symplectic form}) when $\alpha_{0},\alpha_{1}>0$. In this situation, $\omega_{\alpha}$ becomes a K\"ahler form. \\
The extended gauge group  $\widetilde{\mathcal{G}}$ has a canonical action on $\mathcal{A}\times\Omega^{n}(L)$ and this action is symplectic for any $\omega_{\alpha}$. The action is given by 
\[g.(A, \beta)=\left(g.A,\  \widetilde{p(g)}\circ \beta\right),\] 
where $\widetilde{p(g)}$ is the lift of $p(g)\in \mathcal{H}$ to $\widetilde{p(g)}\in \widetilde{\mathcal{H}}$ using the connection on $L$. Using (\ref{Moment map in coupling method}) and (\ref{moment map for KE equation}) , we get the following theorem which provides the moment map interpretation of our equations.
\begin{theorem}
	The $\tilde{\mathcal{G}}$-action on $\mathcal{A}\times\Omega^{n}(L)$ is hamiltonian with respect to $\omega_{\alpha}$ and the $\tilde{\mathcal{G}}$-equivariant moment map $\mu_{\alpha}: \mathcal{A}\times\mathcal{J}_{int}\rightarrow(Lie\tilde{\mathcal{G}})^{*}$ is given by
	 \begin{equation}
	 	\label{13}
	 	\begin{split}	
	 	\langle\mu_{\alpha}(A, [\beta]),\xi\rangle&=4\sqrt{-1}\alpha_{1}\int_{X}tr\left(\theta_{A}\xi\wedge(\sqrt{-1}{\bigwedge}_{\omega} F_{A}-\lambda Id)\right)vol_{\omega}-\int_{X}\phi\frac{\alpha_{0}}{2}\left(\frac{\Omega_{h}}{\int_{X}\Omega_{h}}-\frac{\omega^{n}}{\int_{X}\omega^{n}}\right)\\
	 	&+\int_{X}\phi\left(\alpha_{1}\frac{2}{(n-2)!}tr(F_{A}\wedge F_{A})\wedge\omega^{n-2}+4\alpha_{1}\lambda\frac{1}{(n-1)!}tr(\sqrt{-1}F_{A})\wedge\omega^{n-1}\right), 
	 \end{split}
	 \end{equation}
 where $\xi \in Lie\tilde{\mathcal{G}}$ covers $\eta_{\phi}\in Lie\mathcal{H}$.
\end{theorem}

\section{Obstructions for the existence of solutions}
In this section, we introduce two obstructions for the existence of solutions of the coupled equations (\ref{Coupled equations}). The first one is a Futaki \cite{Futaki} type invariant which relates the existence of a solution to the vanishing of a Lie algebra character. The second one is a Matsushima-Lichnerowicz(\cite{Matsushima}, \cite{Lich 1}) type obstruction which relates the existance of a solutions to the reductivity of the Lie algebra of the automorphism group of the associated principal bundle. 
\subsection{Futaki type invariant as an obstruction}
\label{Futaki invariant section}
We find a Futaki type obstruction to the existence of
solutions of the coupled equations (\ref{Coupled equations}), which generalises the classical Futaki character \cite{Futaki}, following the approach of \cite{Garcia et al}. The moment map interpretation of K\"ahler-Einstein equation in the Fano case \cite{Donaldson 3} gives us an invariant. In the following lemma, we prove that this invariant is the same as the classical Futaki invariant.
\begin{lemma}
	\label{Futaki invariant same }
	Suppose $(X,\omega)$ is a fano manifold and $f$ is the Ricci potential i.e., $Ric(\omega)-\omega=\sqrt{-1}\partial\bar{\partial} f$ normalised so that $\int_{X}^{}e^{f}\omega^{n}=\int_{X}^{}\omega^{n}=C$. Then 
	\begin{equation}
		\int_{X} v(f)\omega^{n}=	\int_{X}\theta_{v}S(\omega)\omega^{n}=C\left(\int_{X}\theta_{v}\left(\frac{\Omega_{h}}{\int_{X}\Omega_{h}}-\frac{\omega^{n}}{C}\right)\right),
	\end{equation} 
	where $\theta_{v}$ is the Hamiltonian function corresponding to a holomorphic vector field $v$(i.e., $i_{v}\omega=\bar{\partial}\theta_{v}$) and $h$ is a Hermitian metric on $K_{X}^{-1}$ whose curvature is $\omega$. We normalize $\theta_{v}$ by requiring $\int_{X}\theta_{v}\omega^{n}=0$.
\end{lemma} 
\begin{proof}
	\begin{equation}	
		\begin{split}
			&\int_{X}\theta_{v}S(\omega)\omega^{n}\\
			&=\int_{X}\theta_{v}n\omega^{n}+\int_{X}n\theta_{v}\sqrt{-1}\partial\bar{\partial}f\wedge\omega^{n-1}\\
			&=\int_{X}n \bar{\partial}\theta_{v}\wedge \sqrt{-1}\partial f\wedge\omega^{n-1}\\
			&=\int_{X} i_{v}\omega^{n}\wedge\sqrt{-1}\partial f\\
			&=\int_{X}v(f)\omega^{n}
		\end{split}
	\end{equation}
	It is a fact that 
	\[\Delta \theta_{v}+\theta_{v}+v(f)=const.\]
	By integrating both sides the constant is seen to be 
	\[const.=\frac{1}{C}\int_{X}\theta_{v}e^{f}\omega^{n}.\]
	Because we have $\int_{X}e^{f}\Delta\theta_{v}\omega^{n}+\int_{X}e^{f}v(f)\omega^{n}=0$. Indeed
	\begin{equation}
		\begin{split}
			&\int_{X}e^{f}\Delta\theta_{v}\omega^{n}\\
			&=-\int_{X}e^{f}\sqrt{-1}\partial f\wedge n \bar{\partial}\theta_{v}\wedge\omega^{n-1}\\
			&=-\int_{X}e^{f}\sqrt{-1}\partial f\wedge i_{v}(\omega)^{n}\\
			&=-\int_{X}e^{f} v(f) \omega^{n}.\\
		\end{split}
	\end{equation}
	So now
	\begin{equation}
		\begin{split}
			&\int_{X}v(f)\omega^{n}\\
			&=\int_{X}\theta_{v}e^{f}\omega^{n}-\int_{X}\theta_{v}\omega^{n}.
		\end{split}
	\end{equation} 
	Now suppose $\omega$ is the curvature of a Hermitian metric $h$. Then it is easy to see that upto a constant $c$,
	\[e^{f}\omega^{n}=c \Omega_{h}.\]
	Again by integration the constant is seen to be 
	\[c=\frac{C}{\int_{X}\Omega_{h}}.\]
	Hence
	\begin{equation}
		\begin{split}
			&\int_{X}v(f)\omega^{n}\\
			&=C\left(\frac{\int_{X}\theta_{v}\Omega_{h}}{\int_{X}\Omega_{h}}-\frac{\int_{X}\theta_{v}\omega^{n}}{C}\right)\\
			&=C\left(\int_{X}\theta_{v}\left(\frac{\Omega_{h}}{\int_{X}\Omega_{h}}-\frac{\omega^{n}}{C}\right)\right).
		\end{split}
	\end{equation}
\end{proof}
Suppose $(E,H)$ is a Hermitian vector bundle of rank $r$ and $\mathcal{E}$ is the associated principal $GL(r,\mathbb{C})$-bundle. Suppose  $I$ is an integrable almost complex structure on $\mathcal{E}$ that preserves the vertical bundle $V\mathcal{E}$ and $\widecheck{I}$ is the unique integrable almost complex structure on the manifold $X$ determined by $I$. Any vector field $\xi \in Lie Aut(\mathcal{E},I)$ covers a real holomorphic vector field $\widecheck{\xi}$ on $(X,\widecheck{I})$ which can be written as 
\begin{equation}
	\widecheck{\xi}=\eta_{\phi_{1}}+\widecheck{I}\eta_{\phi_{2}}+\beta,
\end{equation}
for any K\"ahler form $\omega$ in some K\"ahler class in $(X,\widecheck{I})$, where $\eta_{\phi_{i}}$ is the Hamiltonian vector field of $\phi_{i}$ on $(X,\omega)$ for $i=1,2$, and $\beta$ is the dual of a $1$-form which is harmonic with respect to the K\"ahler metric $\omega(.,\widecheck{I}.)$( see \cite{Lebrun-Simanca}). Since our manifold is Fano, $\beta=0$ in the above expression.  Now we define the Futaki type invariant for our coupled equations, which is the following.
\begin{equation}
	\begin{split}
		&\mathcal{F}_{I}:Lie Aut(\mathcal{E}, I)\rightarrow \mathbb{C}\\
		&\langle\mathcal{F}_{I},\xi\rangle = -4\sqrt{-1}\alpha_{1}\int_{X}tr\left(\theta_{H}\xi\wedge(\sqrt{-1}{\bigwedge}_{\omega}F_{H}-\lambda Id)\right)\frac{\omega^{n}}{n!}\\
		&+\int_{X}\phi \left(\frac{\alpha_{0}}{2}\left(\frac{\Omega_{h}}{\int_{X}\Omega_{h}}-\frac{\omega^{n}}{C}\right)-\alpha_{1}\frac{2}{(n-2)!}tr(F_{H}\wedge F_{H}\wedge\omega^{n-2})-4\alpha_{1}\lambda\frac{1}{(n-1)!}tr(\sqrt{-1}F_{H}\wedge\omega^{n-1})\right),
	\end{split}
\end{equation}
where $\phi=\phi_{1}+\sqrt{-1}\phi_{2}$. One can show that $\mathcal{F}_{I}$ is a character of $Lie Aut(\mathcal{E}, I)$ and does not depend on the choice of $\omega$ and $H$. The proof is essentially same as in \cite{Garcia et al}. We see that $\mathcal{F}_{I}$ vanishes whenever the coupled equations (\ref{Coupled equations}) have a solution.\\
\subsection{Matsushima-Lichnerowicz type obstruction}
\label{matsushima obstruction section}
The first known obstruction for the existence of K\"ahler-Einstein metric on Fano manifold was noticed by Matsushima \cite{Matsushima}. In this section, we introduce the analogous obstruction for the coupled equations (\ref{Coupled equations})
using its moment map interpretaion (\ref{Coupled moment map section}).\par
 We use notations of the previous section(\ref{Futaki invariant section}). Since our manifold $X$ is Fano, we have $H^{1}(X,\mathbb{R})=0$. First we have the following lemma.
\begin{lemma}
	\label{general lemma}
	For any $y\in Lie Aut(\mathcal{E}, I)$, there exist $\xi_{1},\xi_{2}\in Lie\tilde{\mathcal{G}}$ such that $y=\xi_{1}+I\xi_{2}$.
\end{lemma} 
The proof is the same as in (\cite{Garcia et al 1}, lemma $3.7$). 
 \begin{lemma}
 	\label{Lie algebra lemma}
	Suppose $(h,H)$ solves the coupled equations (\ref{Coupled equations}) for positive $\alpha_{0}, \alpha_{1}$. Then for any $y\in Lie Aut(\mathcal{E}, I)$, we have $\xi_{1},\xi_{2}\in Lie Aut(\mathcal{E}, I)$, where $\xi_{1},\xi_{2}$ as in lemma (\ref{general lemma}). 
\end{lemma} 
\begin{proof}
	Since $(h,H)$ is a solution of (\ref{Coupled equations}), we see from section (\ref{Coupled moment map section}) that $\mu_{\alpha}(A,\beta)=0$, for the pair $l=(A,\beta)$. We have seen in section (\ref{Coupled moment map section}) that $\mathcal{A}^{1,1}\times\mathcal{J}_{int}$ is endowed with a formally integrable almost complex structure $J$ and the corresponding K\"ahler form is $\omega_{\alpha}$(we have $\alpha_{0},\alpha_{1}>0$ in our assumption). Now for $y \in  Lie Aut(\mathcal{E}, I)$, we denote by $Y_{y|l}$ the infinitesimal action of $y$ on $l$. First, we see that $Y_{y|l}=0$ since $I$ is integrable. 
	This gives us 
	\begin{equation}
		\label{norm equality}
		\begin{split}
			&0=\lvert\lvert Y_{y|l}\rvert\rvert^{2}\\
			&=\lvert\lvert Y_{\xi_{1}|l}+JY_{\xi_{2}|l}\rvert\rvert^{2}\\
			&=\lvert\lvert Y_{\xi_{1}|l}\rvert\rvert^{2}+\lvert\lvert JY_{\xi_{2}|l}\rvert\rvert^{2}+\langle Y_{\xi_{1}|l},J Y_{\xi_{2}|l}\rangle_{g_{\alpha}}+\langle J Y_{\xi_{2}|l}, Y_{\xi_{1}|l}\rangle_{g_{\alpha}}\\
			&=\lvert\lvert Y_{\xi_{1}|l}\rvert\rvert^{2}+\lvert\lvert Y_{\xi_{2}|l}\rvert\rvert^{2}-2\omega_{\alpha}( Y_{\xi_{1}|l}, Y_{\xi_{2}|l}),
		\end{split}
	\end{equation}
where $g_{\alpha}=\omega_{\alpha}(.,J.)$ is the K\"ahler metric and   $\lvert\lvert.\rvert\rvert$ is the norm on $T_{l}(\mathcal{A}^{1,1}\times \mathcal{J}_{int})$ induced by the metric $g_{\alpha}$. Now
\begin{equation}
	\begin{split}
		&\omega_{\alpha}( Y_{\xi_{1}|l}, Y_{\xi_{2}|l})\\
		&=d\langle \mu_{\alpha},\xi_{1}\rangle(Y_{\xi_{2}|l})\\
		&=\langle\mu_{\alpha}(l),[\xi_{1},\xi_{2}]\rangle\\
		&=0
	\end{split}
\end{equation}
Hence from (\ref{norm equality}), we get that $Y_{\xi_{1}|l}=0=Y_{\xi_{2}|l}$. This implies that $\xi_{1},\xi_{2}\in Lie Aut(\mathcal{E}, I)$. The proof is completed. 
\end{proof}
Now suppose $\tilde{\mathcal{G}}_{l}$ is the isotropy group of $l$ for the action of $\tilde{\mathcal{G}}$ on $\mathcal{A}^{1,1}\times \mathcal{J}_{int}$ and $\mathfrak{k}$ is its Lie algebra. Then we have that 
\[\mathfrak{k}+I\mathfrak{k}\subset Lie Aut(\mathcal{E},I).\]
The Lie group $\tilde{\mathcal{G}}_{l}$ is compact because it can be seen as a closed subgroup of the isometry group of a Riemannian metric on $\mathcal{E}$. 
The lemma (\ref{Lie algebra lemma}) now gives us that 
\[ Lie Aut(\mathcal{E},I)=\mathfrak{k}+I\mathfrak{k}.\]
So $ Lie Aut(\mathcal{E},I)$ is the complexification of the Lie algebra of a compact Lie group. Hence we have the following theorem.
\begin{theorem}
	If $(X,E)$ admits a solution of the coupled equation (\ref{Coupled equations}) for $\alpha_{0},\alpha_{1}>0$, then  $Lie Aut(\mathcal{E},I)$ is reductive.
\end{theorem}
\section{Deformation of solutions}
\label{Deformation section}
In this section, we find solutions of the coupled equations under deformation of the coupling constants. We also find sufficient conditions for the existence of the solutions of the coupled equations. We fix a Fano manifold $X$ and a Hermitian holomorphic vector bundle $E$ over it. We will obtain the solutions by deforming solutions when $\alpha_{0}\neq 0, \alpha_{1}=0$. In this situation, we divide the second equation in (\ref{Coupled equations}) by $\alpha_{0}$ and call $\frac{\alpha_{1}}{\alpha_{0}}$ by $\alpha$.\par
 Let $\mathcal{S}$ be the set of smooth functions $\phi\in C^{\infty}_{0}(X)$ (i.e., $\int_{X}\phi \omega^{n}=0$)  such that $\omega_{\phi}$ is a K\"ahler metric and $\mathcal{W}$ be the set of $Lie{\mathcal{G}}$ elements which are traceless.
 We compute the deformation under the deformations given by \[\omega_{\phi}=\omega+\sqrt{-1}\partial\bar{\partial}\phi\] and \[A_{\xi}=e^{\sqrt{-1}\xi}. A.\]
  Now consider the operator 
\begin{equation}
	\label{moment map operator}
	\begin{split}	
		&T_{\alpha}: \mathcal{S}\times \mathcal{W}\rightarrow C^{\infty}(X)\times Lie{\mathcal{G}}\\
		&\ \ \ \ \ \ (\phi,\xi)\mapsto\left(T^{1}_{\alpha}(\phi,\xi),T^{2}_{\alpha}(\phi,\xi)\right),\\
		&where \\ &T^{1}_{\alpha}(\phi,\xi)=\frac{\frac{1}{2}\left(\frac{\Omega_{h_{\phi}}}{\int_{X}\Omega_{h_{\phi}}}-\frac{\omega_{\phi}^{n}}{\int_{X}\omega_{\phi}^{n}}\right)-\alpha\frac{2}{(n-2)!}tr(F_{A_{\xi}}\wedge F_{A_{\xi}})\wedge\omega_{\phi}^{n-2}-\alpha\tilde{C}vol_{\omega_{\phi}}}{{\omega}^{n}_{\phi}}\\
		&T^{2}_{\alpha}(\phi,\xi)=\sqrt{-1}{\bigwedge}_{\omega_{\phi}}F_{A_{\xi}}-\lambda.
	\end{split}
\end{equation}
Here $h_{\phi}$ is the metric on $K_{X}^{-1}$ corresponding to $\omega_{\phi}$ and $F_{A_{\xi}}$ is the curvature of the connection $A_{\xi}$.\\
Suppose $\omega$ is a K\"ahler-Einstein metric and $A$ is a Hermitian-Yang-Mills connection with respect to $\omega$. The linearisation of this operator for $\alpha=0$ at the point $(\phi=0,\xi=0)$ in the direction $(\dot\phi,\dot\xi)$ is given by
\begin{equation}
	\label{Linearisation of moment map operator}
	\begin{split}
	&\delta T^{1}_{0}(\dot\phi,\dot\xi)=\frac{\frac{1}{2}\left(\frac{-\dot\phi\Omega_{h}}{\int_{X}\Omega_{h}}+\frac{\Omega_{h}\int_{X}\dot\phi\Omega_{h}}{(\int_{X}\Omega_{h})^{2}}-\frac{n{\omega}^{n-1}\wedge\sqrt{-1}\partial\bar{\partial}\dot\phi}{\int_{X}{\omega}^{n}}\right)}{{\omega}^{n}}\\
	&\delta T^{2}_{0}(\dot\phi,\dot\xi)=\sqrt{-1}d^{*}_{A}d_{A}\dot\xi+n(n-1)\frac{\sqrt{-1}F_{A}\wedge\omega^{n-2}\wedge\sqrt{-1}\partial\bar{\partial}\dot\phi}{\omega^{n}}-\frac{n\sqrt{-1}F_{A}\wedge\omega^{n-1}}{\omega^{n}}\times\frac{n\omega^{n-1}\wedge\sqrt{-1}\partial\bar{\partial}\dot\phi}{\omega^{n}}.	
	\end{split}
\end{equation}
 We show that in this situation, the linearisation is an isomorphism under certain assumption. Since the linearisation is elliptic of second order, it is Fredholm. Thus we only need to show that the kernel and cokernel are trivial.\\
 Let us recall that the K\"ahler-Einstein condition is that $\Omega_{h}=\omega^{n}$ \cite{Donaldson 3}. So the first equation in (\ref{Linearisation of moment map operator}) becomes \[\frac{-\Delta_{\omega}\dot\phi-\dot\phi}{2\int_{X}\omega^{n}}.\]
 This operator is invertible if we assume that the manifold has no holomorphic vector field \cite{Lich 2}.  
Let $(\dot\phi,\dot\xi)\in Ker(\delta T_{0})$. Then the previous discussion gives us that $\dot\phi$ is zero. Hence from the second equation in (\ref{Linearisation of moment map operator}) , we get
\[\sqrt{-1}d_{A}^{*}d_{A}\dot\xi=0.\] The $L^{2}$-inner product on $Lie \mathcal{G}$ is given by
\[\langle\xi_{1},\xi_{2}\rangle\coloneqq \int_{X}tr(\xi_{1}\wedge\xi_{2})\frac{\omega^{n}}{n!}.\]
So using this pairing, we have that $d_{A}\dot\xi=0$. Since $\dot\xi\in \mathcal{W}$, we have $\dot\xi=0$. This gives us that the kernel of the linearisation is trivial. \\
Now we show that the cokernel is trivial. To show that we use the following fact. 
\begin{fact}
	Let $I\subset\mathbb{R}$ be a connected subset and $V$ be the set of Fredholm operators. Then for any continuous map $F(t):I\rightarrow V$, the index i.e., \[index(F(t))=dim(ker(F(t))-dim(coker(F(t)))\]is constant. 
\end{fact}
Now consider $I=[0,1]$ and
\begin{equation}
\begin{split}
&F_{1}(t)=\frac{\frac{1}{2}\left(\frac{-\dot\phi\Omega_{h}}{\int_{X}\Omega_{h}}+\frac{\Omega_{h}\int_{X}\dot\phi\Omega_{h}}{(\int_{X}\Omega_{h})^{2}}-\frac{n{\omega}^{n-1}\wedge\sqrt{-1}\partial\bar{\partial}\dot\phi}{\int_{X}\omega_{\phi}^{n}}\right)}{{\omega}^{n}}\\
&F_{2}(t)=\sqrt{-1}d^{*}_{A}d_{A}\dot\xi+tn(n-1)\frac{\sqrt{-1}F_{A}\wedge\omega^{n-2}\wedge\sqrt{-1}\partial\bar{\partial}\dot\phi}{\omega^{n}}-t\frac{n\sqrt{-1}F_{A}\wedge\omega^{n-1}}{\omega^{n}}\times\frac{n\omega^{n-1}\wedge\sqrt{-1}\partial\bar{\partial}\dot\phi}{\omega^{n}}
\end{split}
\end{equation}
At $t=0$, the above operator is an elliptic, self-adjoint and its kernel and cokernel are zero. So using the fact, we have that index at $t=1$ is zero, i.e., at $t=1$, the dimension of kernel and cokernel are the same. But we have shown that the kernel is trivial and hence the cokernel is trivial too. Hence the linearisation is an isomorphism. Hence we have the following theorem.  
\begin{theorem}
	Suppose $\omega$ is a K\"ahler-Einstein metric on a Fano manifold $X$ and $A$ is a Hermitian-Yang-Mills connection on a Hermitian holomorphic vector bundle $E$. Further, assume that $X$ has no non-zero holomorphic vector field. Then there exists $\epsilon>0$ such that for $\tilde{\alpha}\in\mathbb{R}$ with $-\epsilon<\tilde{\alpha}<\epsilon$, there is a solution $(\omega_{\phi},A_{\xi})$ to the coupled equation (\ref{Coupled equations}) with coupling constant $(1,\tilde{\alpha})$. 
\end{theorem} 
\section{Examples of solutions on some projective bundles}
\label{Calabi ansatz section}
In this section, we use Calabi ansatz to produce solutions to our equations. We follow the method of Keller-T{\o}nnesen-Friedman \cite{Keller} while the calculations in S\'zekelyhidi \cite{Gabor}.\\
Let $\Sigma_{i}=\mathbb{C}P^{1}$ and consider the line bundle $L=\otimes_{i=1}^{k}\pi_{i}^{*}(L_{i})$ over $\prod_{i=1}^{k}\Sigma_{i}$, where $L_{i}$ is a holomorphic line bundle of degree $-1$ over $\Sigma_{i}$. Set $X_{k}:=\mathbb{P}(L \oplus \mathcal{O})\rightarrow \prod_{i=1}^{k}\Sigma_{i}$, where $\mathcal{O}$ is the trivial line bundle over $\prod_{i=1}^{k}\Sigma_{i}$. We denote $L\oplus \mathcal{O}$ by $E$. So $\det(E)\cong L$.\par
 Let $\omega_{\Sigma_{i}}$ be a K\"ahler metric on $\Sigma_{i}$ with constant scalar curvature 2. By the Gauss-Bonnet theorem, the area of $\Sigma_{i}$ is $2\pi$ with this metric. Let $h'$ be a metric on $L$ with curvature form 
 	$F(h')=-\sum_{i=1}^{k}\omega_{\Sigma_{i}}$. We consider metric of the form
 \begin{equation}
  \omega_{k}=\sum_{i=1}^{k}p_{i}^{*}\omega_{\Sigma_{i}}+\sqrt{-1}\partial\bar{\partial}f_{k}(s),
\end{equation}
  where $p_{i}:L\rightarrow \Sigma_{i}$ is the projection map and $s=\log(\lvert(z_{1},\dots,z_{k},w)\rvert_{h'}^{2})$. We suppress the $k$ in the subsequent calculations. \par
  Now $\lvert(z_{1},\dots,z_{k},w)\rvert_{h'}^{2}=\lvert w\rvert^{2} h'(z_{1},\dots,z_{k})$, so $s=\log\lvert w\rvert^{2}+\log h'(z_{1},\dots,z_{k})$. We work at a point $(z_{1}^{0},\dots,z_{k}^{0},w^{0})$ where $d\log h'(z_{1}^{0},\dots,z_{k}^{0})=0$. At this point, $\sqrt{-1}\partial\bar{\partial}f(s)=f'(s) \big(\sum_{i=1}^{k}p_{i}^{*}\omega_{\Sigma_{i}}\big)+\sqrt{-1}f''(s)\frac{dw\wedge d\bar{w}}{\lvert w\rvert^{2}}$. So
  \begin{equation}
  	\label{Our ansatz}
   \omega=\big(1+f'(s)\big)\big(\sum_{i=1}^{k}p_{i}^{*}\omega_{\Sigma_{i}}\big)+\sqrt{-1}f''(s)\frac{dw\wedge d\bar{w}}{\lvert w\rvert^{2}}.
\end{equation}
We still need to extend the metric across the zero and infinity sections. The metric in the fiber directions is given by $\sqrt{-1}f''(s)\frac{dw\wedge d\bar{w}}{\lvert w\rvert^{2}}$. We know that for our ansatz (\ref{Our ansatz}) to be positive, $f$ must be strictly convex. Hence we can take the Legendre transform of $f$. The Legendre transform $F$ is defined in terms of the variable $\tau=f'(s)$ by the formula
\[f(s)+F(\tau)=s\tau.\]
If the range of $\tau$ is an interval $I\subset \mathbb{R}$, then $F$ is a strictly convex function defined on $I$. The momentum profile of the metric is defined to be $\phi: I \rightarrow \mathbb{R}$, where 
\[\phi(\tau)=\frac{1}{F''(\tau)}.\] 
We have the following useful relations:
\begin{equation}
	\label{Legendre transform relation}
	s=F'(\tau),\ \ \ \  \frac{ds}{d\tau}=F''(\tau),\ \ \ \  \phi(\tau)=f''(s).
\end{equation}
Now we can write (\ref{Our ansatz}) as
\begin{equation}
	\label{Our ansatz in Legendre transform form}
	 \omega=\big(1+\tau\big)\big(\sum_{i=1}^{k}p_{i}^{*}\omega_{\Sigma_{i}}\big)+\sqrt{-1}\phi(\tau)\frac{dw\wedge d\bar{w}}{\lvert w\rvert^{2}}.
\end{equation}
  It turns out that the conditions for the metric to be extended across the zero and infinity sections are-
\begin{equation}
	\label{Boundary conditions}
	\begin{split}
		& \lim_{\tau\rightarrow a}\phi(\tau)=0, \ \ \ \lim_{\tau\rightarrow a}\phi'(\tau)=1\\
		& \lim_{\tau\rightarrow b}\phi(\tau)=0, \ \ \ \lim_{\tau\rightarrow b}\phi'(\tau)=-1,
	\end{split}
\end{equation}
where $(a,b)$ is the range of $\tau$.\par
From (\ref{Our ansatz in Legendre transform form}) , we see that the metric will be positive definite as long as $1+\tau$ and $\phi(\tau)$ are positive on $(a,b)$. For simplicity, we can take the interval $[0,m]$ for some $m>0$. The value of $m$ determines the K\"ahler class of the resulting metric.\par
Now $K_{X_{k}}^{-1}=\mathcal{O}_{X_{k}}(2)\otimes \pi^{*}(det E \otimes K_{Y_{k}}^{-1})$(see lemma $2.1$ in \cite{Zhang}), where $Y_{k}=\prod_{i=1}^{k}\Sigma_{i}$ and $\pi:X_{k}\rightarrow Y_{k}$ the projection map. We want to define a metric $h=h_{0}e^{-f}$ on $K_{X_{k}}^{-1}$ such that $F(h)=\omega=(1+f'(s))\big(\sum_{i=1}^{k}p_{i}^{*}\omega_{\Sigma_{i}}\big)+\sqrt{-1}f''(s)\frac{dw\wedge d\bar{w}}{\lvert w\rvert^{2}}.$ To do that we first need to find the line bundle where $e^{-f}$ is a metric on it. \par
Viewing $X_{k}$ as a $\mathbb{C}P^{1}$ bundle over $\prod_{i=1}^{k}\Sigma_{i}$, we see that the space $H^{2}(X_{k},\mathbb{R})$ is generated by Poincar\'e duals of $D_{\Sigma_{i}}$ and $S_{0}$, where the zero section $S_{0}$ is the image of the subbundle $\{0\}\oplus\mathcal{O}$ under the projection map to the projectivization $X_{k}=\mathbb{P}(L \oplus \mathcal{O})$ and $D_{\Sigma_{i}}=\mathbb{P}(E)|_{\Sigma_{1}\times\dots\times[1:0]\times\dots\times\Sigma_{k}}$. $D_{\Sigma_{i}}$ are divisors of $\mathbb{P}(E)$ corresponding to the line bundles $(\pi_{i}\circ\pi)^{*}(\mathcal{O}_{\Sigma_{i}}(1))$ and the divisor $S_{0}$ corresponds to the line bundle $\mathcal{O}_{X_{k}}(1)\otimes \pi^{*}L$(see lemma $2.1$ in \cite{Zhang}). We suppress the pullbacks from now on. We denote the infinity section by $S_{\infty}$, which is the image of the subbundle $L \oplus\{0\}$ under the projection map to the projectivization $X_{k}=\mathbb{P}(L \oplus \mathcal{O})$.\par
Suppose $a_{0}[S_{0}]+\sum_{i=1}^{k}a_{i}[D_{\Sigma_{i}}]$ be the Poincar\'e
 dual of the K\"ahler class of $\sqrt{-1}\partial\bar{\partial}f(s)$, then the coefficients can be calculated by integration.  Indeed,
 \begin{equation}
 a_{0}=\int_{\cap_{i=1}^{k}D_{\Sigma_{i}}}^{}\sqrt{-1}\partial\bar{\partial}f(s)=\int_{\mathbb{C}\setminus {0}}^{}\sqrt{-1}f''(s)\frac{dw\wedge d\bar{w}}{\lvert w\rvert^{2}}=2\pi\big(\lim_{s\rightarrow\infty}f'(s)-\lim_{s\rightarrow -\infty}f'(s)\big)=2\pi m
\end{equation}
\begin{equation}
	a_{i}=\int_{S_{\infty}\cap_{j=1,j\neq i}^{k} D_{\Sigma_{j}}}^{}\sqrt{-1}\partial\bar{\partial}f(s)=\int_{S_{\infty}\cap_{j=1,j\neq i}^{k} D_{\Sigma_{j}}}^{}f'(s)\omega_{\Sigma_{i}}=2\pi m,
\end{equation}
where we have used our assumptions that $\lim\limits_{s\rightarrow\infty}f'(s)=m$ and $\lim\limits_{s\rightarrow-\infty}f'(s)=0$.\par
Now suppose $m=2$. Since $S_{0}$ corresponds to the line bundle $\mathcal{O}_{X_{k}}(1)\otimes L$, we see that $e^{-f}$ is a metric on the line bundle $(\mathcal{O}_{X_{k}}(1)\otimes L)^{2}\otimes_{i=1}^{k}\mathcal{O}_{\Sigma_{i}}(2)$. 
Now $K_{X_{k}}^{-1}$ can be written as 
\begin{equation}
	K_{X_{k}}^{-1}=\mathcal{O}_{X_{k}}(2)\otimes L^{2}\otimes_{i=1}^{k} \mathcal{O}_{\Sigma_{i}}(2)\otimes L^{*}.
\end{equation}
Suppose $\tilde{h'}$ is the dual metric corresponding to $h'$ on $L^{*}$, then its curvature is $F(\tilde{h'})=\sum_{i=1}^{k}\omega_{\Sigma_{i}}$. Now, it is clear that $h=\tilde{h'}e^{-f}$ is a metric on $K_{X_{k}}^{-1}$ with curvature form $F(h)=\sum_{i=1}^{k}\omega_{\Sigma_{i}}+\sqrt{-1}\partial\bar{\partial}f=\omega$.\par
This metric $h$ can be viewed as a volume form $\Omega_{h}$ which in terms of local coordinates looks like $\Omega_{h}=V(\sqrt{-1})^{k+1}dz^{1}\wedge d\bar{z}^{1}\wedge\dots\wedge dz^{k}\wedge d\bar{z}^{k}\wedge dw\wedge d\bar{w}$, where $V=\lvert\frac{\partial}{\partial z^{1}}\wedge\dots\wedge\frac{\partial}{\partial z^{k}}\wedge\frac{\partial}{\partial w}\rvert_{h}^{2}$. We can see from the isomorphism $K_{X_{k}}^{-1}\cong \mathcal{O}_{X_{k}}(2)\otimes\pi^{*}(\det E\otimes K_{Y_{k}}^{-1})$ that \begin{equation}
	\label{V expression}
V=\tilde{h'} e^{-f}.
\end{equation}
The Fubini-Study metric in local coordinates can be written as 
\begin{equation}
	\label{FS expression}
	\omega_{\Sigma_{i}}=\sqrt{-1}\tilde{\omega}_{\Sigma_{i}}dz^{i}\wedge d\bar{z}^{i}.
\end{equation} 
Now using (\ref{V expression}) and (\ref{FS expression}), we get 
\begin{equation}
	\label{volume expression}
	\Omega_{h}=\frac{\tilde{h'}e^{-f}}{\tilde{\omega}_{\Sigma_{1}}\dots\tilde{\omega_{\Sigma_{k}}}}\omega_{\Sigma_{1}}\wedge\dots\wedge\omega_{\Sigma_{k}}\wedge \sqrt{-1}dw\wedge d\bar{w}
\end{equation}
 Now $\tilde{h'}$ is a metric on $L^{*}\cong \otimes_{i=1}^{k}\mathcal{O}_{\Sigma_{i}}(1)$ and $\omega_{\Sigma_{i}}$ is a metric on $\mathcal{O}_{\Sigma_{i}}(2)$. So $\frac{\tilde{h'}}{\tilde{\omega_{\Sigma_{1}}}\dots\tilde{\omega_{\Sigma_{2}}}}$ is a metric on $L$ and it is the same as $h'$. So now (\ref{volume expression}) becomes 
\begin{equation}
	\label{expression for Omega h}
	\begin{split}
	\Omega_{h}&=h'e^{-f}\omega_{\Sigma_{1}}\wedge\dots\wedge\omega_{\Sigma_{k}}\wedge\sqrt{-1}dw\wedge d\bar{w}\\
	&=\lvert w\rvert^{2}h'e^{-f}\omega_{\Sigma_{1}}\wedge\dots\wedge\omega_{\Sigma_{k}}\wedge\sqrt{-1}\frac{dw\wedge d\bar{w}}{\lvert w\rvert^{2}}\\
	&=e^{s-f}\omega_{\Sigma_{1}}\wedge\dots\wedge\omega_{\Sigma_{k}}\wedge\sqrt{-1}\frac{dw\wedge d\bar{w}}{\lvert w\rvert^{2}}.
	\end{split}
\end{equation}
Now we will calculate the volume form associated with our ansatz. Our ansatz is $\omega=(1+f'(s))\sum_{i=1}^{k}\omega_{\Sigma_{i}}+\sqrt{-1}f''(s)\frac{dw\wedge d\bar{w}}{\lvert w\rvert^{2}}$. 
So
\begin{equation}
	\label{expression for omega to the power k+1}
	\omega^{k+1}=(k+1)!(1+f'(s))^{k}f''(s)\omega_{\Sigma_{1}}\wedge\dots\wedge\omega_{\Sigma_{k}}\wedge\frac{\sqrt{-1}dw\wedge d\bar{w}}{\lvert w\rvert^{2}}
\end{equation}
and hence 
\begin{equation}
	\label{expression of volume form}
	Vol_{\omega}=(1+f'(s))^{k}f''(s)\omega_{\Sigma_{1}}\wedge\dots\wedge\omega_{\Sigma_{k}}\wedge\frac{\sqrt{-1}dw\wedge d\bar{w}}{\lvert w\rvert^{2}}.
\end{equation}
 Now we calculate $\omega^{k-1}$.
\begin{equation}
	\label{expression for omega k-1}
	\begin{split}
		\omega^{k-1}&=\left(1+f'(s)\right)^{k-1}\left(\sum_{i=1}^{k}\omega_{\Sigma_{i}}\right)^{k-1}+(k-1)\left(1+f'(s)\right)^{k-2}\left(\sum_{i=1}^{k}\omega_{\Sigma_{i}}\right)^{k-2}f''(s)\frac{\sqrt{-1}dw\wedge d\bar{w}}{\lvert w\rvert^{2}}\\
		&=\left(1+f'(s)\right)^{k-1}(k-1)!\sum_{i=1}^{k}\omega_{\Sigma_{1}}\wedge\dots\wedge\widehat{\omega_{\Sigma_{i}}}\wedge\dots\wedge\omega_{\Sigma_{k}}\\
		&+(k-1)!\left(1+f'(s)\right)^{k-2}f''(s)\frac{\sqrt{-1}dw\wedge d\bar{w}}{\lvert w\rvert^{2}}\sum_{i<j}^{}\omega_{\Sigma_{1}}\wedge\dots\wedge\widehat{\omega_{\Sigma_{i}}}\wedge\dots\wedge\widehat{\omega_{\Sigma_{j}}}\wedge\dots\wedge\omega_{\Sigma_{k}}.
	\end{split}
\end{equation}

Now we consider two cases.
\subsection{Case 1: Closed, traceless $(1,1)$ form of general type}

In this case, we take $\gamma_{k}=\sum_{i=1}^{k}\omega_{\Sigma_{i}}-k\frac{\phi(\tau)}{1+\tau}\frac{\sqrt{-1}dw\wedge d\bar{w}}{\lvert w\rvert^{2}}$ as our closed, traceless $(1,1)$ form. $\gamma_{k}$ is globally defined because $\tau$ and $\phi(\tau)$ are globally defined functions. For integers $(m_{1},m_{2})\neq(0,0)$, we take $a_{k}=\frac{m_{1}+km_{2}\log3}{2+3k\log3}$ and $b_{k}=\frac{2m_{2}-3m_{1}}{2+3k\log3}$. Then $\frac{a_{k}\omega+b_{k}\gamma_{k}}{2\pi}$ is an integral cohomology class and hence there exists a Hermitian line bundle $(L_{k},H_{k})$ such that $\sqrt{-1}F_{H_{k}}=\frac{a_{k}\omega+b_{k}\gamma_{k}}{2\pi}$.\par  
 This gives $\sqrt{-1}\bigwedge_{\omega}F_{H_{k}}=\frac{(k+1)a_{k}}{2\pi}$ and hence the first equation of our coupled equation (\ref{Coupled equations})  is met with $\lambda=\frac{(k+1)a_{k}}{2\pi}$.\par
Now for the second equation, we need to calculate the term $F_{H_{k}}\wedge F_{H_{k}}\wedge\omega^{k-1}$. In this case
\begin{equation}
	\label{Second Chern character for general alpha}
	\begin{split}
		&F_{H_{k}}\wedge F_{H_{k}}\wedge\omega^{k-1}\\
		&=-\frac{a_{k}^{2}\omega^{k+1}+b_{k}^{2}\gamma_{k}^{2}\wedge\omega^{k-1}}{(2\pi)^{2}}\\
	\end{split}
\end{equation} 
We see that 
\begin{equation}
	\label{second term}
	\gamma_{k}^{2}=(\sum_{i=1}^{k}\omega_{\Sigma_{i}})^{2}-2k\frac{\phi(\tau)}{1+\tau}(\sum_{i=1}^{k}\omega_{\Sigma_{i}})\wedge\frac{\sqrt{-1}dw\wedge d\bar{w}}{\lvert w\rvert^{2}}.
\end{equation}
Now putting (\ref{expression for omega to the power k+1}), (\ref{expression for omega k-1}) and  (\ref{second term}) in (\ref{Second Chern character for general alpha}), we get
\begin{equation}
	\label{Second Chern character for general alpha final form}
	\begin{split}
		&F_{H_{k}}\wedge F_{H_{k}}\wedge\omega^{k-1}\\
		&=-\frac{a_{k}^{2}(k+1)!(1+\tau)^{k}-b_{k}^{2}(k+1)!(1+\tau)^{k-2}}{(2\pi)^{2}}\phi(\tau)\omega_{\Sigma_{1}}\wedge\dots\wedge\omega_{\Sigma_{k}}\wedge\frac{\sqrt{-1}dw\wedge d\bar{w}}{\lvert w\rvert^{2}}
	\end{split}
\end{equation}
Putting  (\ref{expression for Omega h}), (\ref{expression for omega to the power k+1}),   (\ref{expression of volume form}), (\ref{Second Chern character for general alpha final form}) and $\int_{X_{k}}\Omega_{h}=C_{k}'$, $\int_{X_{k}}\omega^{k+1}=C_{k}$ in the second equation of (\ref{Coupled equations}) , we get 
\begin{equation}
	\label{k less or equal to 4 case equation}
	\begin{split}
		&\frac{\alpha_{0}}{2}\left(\frac{e^{s-f(s)}}{C_{k}'}-\frac{(k+1)!(1+\tau)^{k}\phi(\tau)}{C_{k}}\right)+\frac{2\alpha_{1}}{(k-1)!(2\pi)^{2}}\left(a_{k}^{2}(k+1)!(1+\tau)^{k}-b_{k}^{2}(k+1)!(1+\tau)^{k-2}\right)\phi(\tau)\\
		&=\tilde{C_{k}}(1+\tau)^{k}\phi(\tau)\\
		&\implies\frac{\alpha_{0}e^{s-f(s)}}{2C_{k}'}=\phi(\tau)\left[\left(\frac{\alpha_{0}(k+1)!}{2C_{k}}-\frac{2\alpha_{1}(k+1)ka_{k}^{2}}{(2\pi)^{2}}+\tilde{C_{k}}\right)(1+\tau)^{k}+\frac{b_{k}^{2}2\alpha_{1}(k+1)k}{(2\pi)^{2}}(1+\tau)^{k-2}\right].
	\end{split}
\end{equation}
Now we set 
\begin{equation}
	\label{G k tau }
	G_{k}(\tau)=\left(\frac{\alpha_{0}(k+1)!}{2C_{k}}-\frac{2\alpha_{1}(k+1)ka_{k}^{2}}{(2\pi)^{2}}+\tilde{C_{k}}\right)(1+\tau)^{k}+\frac{b_{k}^{2}2\alpha_{1}(k+1)k}{(2\pi)^{2}}(1+\tau)^{k-2}.
\end{equation}
Assume that $\alpha_{0}>0$ and $G_{k}(\tau)>0$ for $\tau\in [0,2]$. Then taking $\log$ on both side of (\ref{k less or equal to 4 case equation}) and differentiating with respect to $\tau$ we get 
\begin{equation}
	\begin{split}
		&(1-\tau)G_{k}(\tau)=\left[\phi(\tau)G_{k}(\tau)\right]'\\
		&\phi(\tau)=\frac{\int_{0}^{\tau}(1-t)G_{k}(t)dt}{G_{k}(\tau)}
	\end{split}
\end{equation}
In the first line, we have used (\ref{Legendre transform relation}). It is easy to see that $\phi(0)=0$ and $\phi'(0)=1$. Now if we can show that $\phi(2)=0$, then it also implies $\phi'(2)=-1$. Now $\phi(2)=0$ implies
 \begin{equation}
 	\label{phi(2)=0 condition}
 	\begin{split}
 		&\int_{0}^{2}(1-t)G_{k}(t)dt=0\\
 		&\implies\frac{\alpha_{0}(k+1)!}{2C_{k}}-\frac{2\alpha_{1}(k+1)ka_{k}^{2}}{(2\pi)^{2}}+\tilde{C_{k}}=-\frac{b_{k}^{2}2\alpha_{1}(k+1)k}{(2\pi)^{2}}\frac{\int_{0}^{2}(1-t)(1+t)^{k-2}dt}{\int_{0}^{2}(1-t)(1+t)^{k}dt}.
 	\end{split}
 \end{equation}
We set 
\begin{equation}
	\label{ R k}
	R(k)=\frac{\int_{0}^{2}(1-t)(1+t)^{k-2}dt}{\int_{0}^{2}(1-t)(1+t)^{k}dt}.
\end{equation}
Now putting (\ref{phi(2)=0 condition}) in (\ref{G k tau }), we get
\begin{equation}
	\label{positivity of G k tau}
	G_{k}(\tau)=\frac{2\alpha_{1}b_{k}^{2}(k+1)k}{(2\pi)^{2}}(1+\tau)^{k-2}\left(1-(1+\tau)^{2}R_{k}\right)
\end{equation}
\subsubsection{Subcase 1: $k=1,2,3,4$}
Simple calculations gives us 
\begin{equation}
	\label{R k calculations}
	R_{1}=-3(\log3-1),\  R_{2}=0,\ R_{3}=\frac{5}{63},\  R_{4}=\frac{5}{46}.
\end{equation}
Now from (\ref{positivity of G k tau}), we see that for $k=1,2,3,4$; $G_{k}(\tau)$ is positive on $[0,2]$ whenever $\alpha_{1}$ is positive . This means that we have solutions in these cases but we need to calculate $C_{k}$ and $\tilde{C_{k}}$ to make the condition (\ref{phi(2)=0 condition}) more compact.\par
We now calculate $C_{k}$ and $\tilde{C_{k}}$.  Writing the fiber coordinate as $w=re^{\sqrt{-1}\theta}$, we have $\frac{\sqrt{-1}dw\wedge d\bar{w}}{\lvert w\rvert^{2}}=2\frac{dr}{r}\wedge d\theta$. From the relation $s=\log (\lvert w\rvert^{2})$, we have $ds=2\frac{dr}{r}$. From (\ref{Legendre transform relation}) , we have $2\frac{dr}{r}=ds=\frac{d\tau}{\phi(\tau)}$. Hence we get the important relation
\begin{equation}
	\label{Fiber volume form relation}
	\frac{\sqrt{-1}dw\wedge d\bar{w}}{\lvert w\rvert^{2}}=ds\wedge d\theta=\frac{d\tau}{\phi(\tau)}\wedge d\theta.
\end{equation}
We use the relations (\ref{Fiber volume form relation}) to calculate $C_{k}$ and $\tilde{C_{k}}$.
\begin{equation}
	\label{C_{k} expression}
	\begin{split}
		C_{k}&=\int_{X_{k}}\omega^{k+1}\\
		&=\int_{X_{k}}(k+1)!(1+\tau)^{k}\phi(\tau)\omega_{\Sigma_{1}}\wedge\dots\wedge\omega_{\Sigma_{k}}\wedge\frac{\sqrt{-1}dw\wedge d\bar{w}}{\lvert w\rvert^{2}}\\
		&=(k+1)!(2\pi)^{k+1}\frac{3^{k+1}-1}{k+1}
	\end{split}
\end{equation}
and 
	\[\tilde{C_{k}}=-\frac{2\alpha_{1}}{(k-1)!}\frac{\int_{X_{k}}F_{H}\wedge F_{H}\wedge\omega^{k-1}}{\frac{C_{k}}{(k+1)!}}\]
Using (\ref{Second Chern character for general alpha final form}) and (\ref{C_{k} expression}) , we get 
\begin{equation}
	\begin{split}
		&\tilde{C_{k}}\\
		&=\frac{2\alpha_{1}(k+1)}{(k-1)!(2\pi)^{k+3}(3^{k+1}-1)}\int_{X_{k}}\left(a_{k}^{2}(k+1)!(1+\tau)^{k}-b_{k}^{2}(k+1)!(1+\tau)^{k-2}\right)\phi(\tau)\bigwedge_{i=1}^{k}\omega_{\Sigma_{i}}\wedge\frac{\sqrt{-1}dw\wedge d\bar{w}}{\lvert w\rvert^{2}}\\
		&=\frac{2\alpha_{1}(k+1)^{2}k}{(2\pi)^{2}(3^{k+1}-1)}\int_{0}^{2}\left(a_{k}^{2}(1+\tau)^{k}-b_{k}^{2}(1+\tau)^{k-2}\right)d\tau.
	\end{split}
	\end{equation}
Hence we have $\tilde{C_{1}}=\frac{\alpha_{1}}{(2\pi)^{2}}(4a_{k}^{2}-b_{k}^{2}\log3)$ and $\tilde{C_{k}}=\frac{2\alpha_{1}(k+1)^{2}k}{(2\pi)^{2}(3^{k+1}-1)}\left(a_{k}^{2}\frac{3^{k+1}-1}{k+1}-b_{k}^{2}\frac{3^{k-1}-1}{k-1}\right)$ for $k\geq 2$.
Hence the condition (\ref{phi(2)=0 condition}) for $k=1$ is \[\alpha_{0}=8b_{1}^{2}\alpha_{1}(13\log3-12)\] 
and for $k=2,3,4$ is 
\[\frac{k+1}{2(2\pi)^{k+1}(3^{k+1}-1)}\alpha_{0}=\alpha_{1}\frac{2b_{k}^{2}(k+1)k}{(2\pi)^{2}}\left(\frac{(k+1)(3^{k-1}-1)}{(k-1)(3^{k+1}-1)}-R_{k}\right)\]
\begin{theorem}
	Suppose $\phi_{k}:[0,2]\rightarrow \mathbb{R}$ is the function $\phi_{k}(x)=\frac{\int_{0}^{x}(1-t)G_{k}(t)dt}{G_{k}(x)}$, where $G_{k}(x)$ is given by (\ref{G k tau }). Then by construction, there is a solution of the coupled equations satisfying the following conditions-
	\begin{equation}
		\begin{split}
			&1) \alpha_{0},\alpha_{1}>0\\
			&2) \alpha_{0}=8b_{1}^{2}\alpha_{1}(13\log3-12)\ \ (for\  k=1)\\
			&3)\frac{k+1}{2(2\pi)^{k+1}(3^{k+1}-1)}\alpha_{0}=\alpha_{1}\frac{2b_{k}^{2}(k+1)k}{(2\pi)^{2}}\left(\frac{(k+1)(3^{k-1}-1)}{(k-1)(3^{k+1}-1)}-R_{k}\right)  (for\  k=2,3,4)
		\end{split}
	\end{equation}
\end{theorem}
\subsubsection{Subcase 2: $k\geq 5$}
Now we see what happens when $k\geq 5$. In this case, we have the following lemma.
\begin{lemma}
	\label{Lemma}
	$1>R(k)=\frac{\int_{0}^{2}(1-t)(1+t)^{k-2}dt}{\int_{0}^{2}(1-t)(1+t)^{k}dt}>\frac{1}{9}$ for $k\geq 5$ and it goes to $\frac{1}{9}$ as $k$ goes to infinity.
\end{lemma}
\begin{proof}
	\begin{equation}
		\label{A_{k}}
		\begin{split}
			&A_{k}=\int_{0}^{2}(1-t)(1+t)^{k}dt\\
			&=\int_{1}^{3}(2-s)s^{k}ds\\
			&=\frac{2}{k+1}(3^{k+1}-1)+\frac{1}{k+2}(1-3^{k+2})\\
			&=\frac{3^{k+1}(1-k)}{(k+1)(k+2)}-\frac{k+3}{(k+1)(k+2)}.
		\end{split}
	\end{equation}
	It is easy to see that $A_{k}<0$ for $k\geq 5$. Now
	\begin{equation}
		\begin{split}
			&A_{k}-A_{k-2}\\
			&=3^{k-1}\times8k[\frac{-k^{2}+4k+3}{(k-1)k(k+1)(k+2)}]-2\times3^{k-1}[\frac{10}{(k-1)(k+1)}+\frac{3}{k(k+2)}]+2[\frac{k^{2}+4k+1}{(k-1)k(k+1)(k+2)}]\\
			&<0
		\end{split}
	\end{equation}
	So we have $1>R_{k}$. Similarly,
	\begin{equation}
		\begin{split}
			&A_{k}-9A_{k-2}\\
			&=2\times3^{k+1}[\frac{k^{2}-4k-3}{(k-1)k(k+1)(k+2)}]+8k[\frac{k^{2}+4k+1}{(k-1)k(k+1)(k+2)}]+[\frac{2k^{2}+40k+18}{(k-1)k(k+1)(k+2)}]\\
			&>0.
		\end{split}
	\end{equation}
	Hence $R_{k}>\frac{1}{9}$.\\
	Using (\ref{A_{k}}) we get
	\begin{equation}
		\begin{split}
			&R_{k}\\
			&=\frac{\frac{2}{k-1}(3^{k-1}-1)+\frac{1}{k}(1-3^{k})}{\frac{2}{k+1}(3^{k+1}-1)+\frac{1}{k+2}(1-3^{k+2})}\\
			&=\frac{1}{9}\times\frac{\frac{2}{k-1}(3^{k-1}-1)+\frac{1}{k}(1-3^{k})}{\frac{2}{k+1}(3^{k-1}-\frac{1}{9})+\frac{1}{k+2}(\frac{1}{9}-3^{k})}.
		\end{split}
	\end{equation}
	From the above expression it is clear that $R_{k}\rightarrow \frac{1}{9}$ as $k\rightarrow\infty$.
\end{proof}
Using lemma (\ref{Lemma}), we see that (\ref{positivity of G k tau}) can not have a definite sign in $[0,2]$ for any choice of $\alpha_{1}$. Hence in this case, our method does not produce any solution.
\subsection{Case 2: Closed, traceless $(1,1)$ form of specific type}
In this case, we choose $\gamma_{k}$ to be a combination of differential forms coming from the base only.
Let $\gamma_{k}$ be a closed, traceless $(1,1)$ form. Consider the line bundle $L_{k}$ with curvature form $\sqrt{-1}F_{H_{k}}=\frac{\omega+\gamma_{k}}{2\pi}$. \par
Then we have 
\begin{equation}
	\label{second Chern character }
	\begin{split}
		&\frac{2}{(k-1)!}F_{H}\wedge F_{H}\wedge \omega^{k-1}\\
		&=-\frac{2}{(2\pi)^{2}(k-1)!}(\omega^{k+1}+\gamma_{k}^{2}\wedge\omega^{k-1}),
	\end{split}
\end{equation}
where we have used the fact that $\gamma_{k}$ is traceless.
\subsubsection{Subcase 1: k even}
 We take $\gamma_{k}=\sum_{i=1}^{k}(-1)^{i-1}\omega_{\Sigma_{i}}$ as our traceless, closed, $(1,1)$ form. One can check that, it is indeed traceless and closed. The cohomology class $\frac{\omega+\gamma_{k}}{2\pi}$ is integral. Hence there exists a line bundle $L_{k}$ with curvature form $\sqrt{-1}F_{H_{k}}=\frac{\omega+\gamma_{k}}{2\pi}.$\par
Now, $k!\sqrt{-1}\bigwedge_{\omega}F_{H_{k}}\frac{\omega^{k+1}}{(k+1)!}=\sqrt{-1}F_{H_{k}}\wedge\omega^{k}=\frac{(\omega+\gamma_{k})}{2\pi}\wedge\omega^{k}=\frac{\omega^{k+1}}{2\pi}$. So $\sqrt{-1}\bigwedge_{\omega}F_{H_{k}}=\frac{k+1}{2\pi}$. Hence, the first equation is met with $\lambda=\frac{k+1}{2\pi}$.\\
Now
\begin{equation}
	\label{expression for alpha square in even case}
	 \gamma_{k}^{2}=2\sum_{i<j}^{}(-1)^{i+j-2}\omega_{\Sigma_{i}}\wedge\omega_{\Sigma_{j}}
\end{equation}
Using (\ref{expression for omega k-1}) and (\ref{expression for alpha square in even case}), we have
\begin{equation}
	\label{alpha square omega to the power k-1}
	\begin{split}
		\gamma_{k}^{2}\wedge\omega^{k-1}&=-k.(k-1)!(1+f'(s))^{k-2}f''(s)\omega_{\Sigma_{1}}\wedge\dots\wedge\omega_{\Sigma_{k}}\wedge\frac{\sqrt{-1}dw\wedge d\bar{w}}{\lvert w\rvert^{2}}\\
		&=-k!(1+f'(s))^{k-2}f''(s)\omega_{\Sigma_{1}}\wedge\dots\wedge\omega_{\Sigma_{k}}\wedge\frac{\sqrt{-1}dw\wedge d\bar{w}}{\lvert w\rvert^{2}}
	\end{split}
\end{equation}
Putting (\ref{expression for omega to the power k+1}) and (\ref{alpha square omega to the power k-1}) in (\ref{second Chern character }), we get
\begin{equation}
	\label{second Chern character in even case}
	\frac{2}{(k-1)!}F_{H}\wedge F_{H}\wedge \omega^{k-1}=-\frac{2}{(2\pi)^{2}(k-1)!}[(k+1)!(1+f'(s))^{k}f''(s)-k!(1+f'(s))^{k-2}f''(s)]\bigwedge_{i=1}^{k}\omega_{\Sigma_{i}}\wedge\frac{\sqrt{-1}dw\wedge d\bar{w}}{\lvert w\rvert^{2}}
\end{equation}
Now putting  (\ref{expression for Omega h}), (\ref{expression for omega to the power k+1}),   (\ref{expression of volume form}), (\ref{second Chern character in even case}) and $\int_{X_{k}}\Omega_{h}=C_{k}'$, $\int_{X_{k}}\omega^{n}=C_{k}$ in the second equation of (\ref{Coupled equations}) , we get 
\begin{equation}
	\label{K even equation 1}
	\begin{split}
&\frac{\alpha_{0}}{2}[\frac{e^{s-f(s)}}{C_{k}'}-\frac{(k+1)!(1+f'(s))^{k}f''(s)}{C_{k}}]+\frac{2\alpha_{1}}{(2\pi)^{2}(k-1)!}[(k+1)!(1+f'(s))^{k}f''(s)-k!(1+f'(s))^{k-2}f''(s)]\\
&=\tilde{C_{k}}(1+f'(s))^{k}f''(s)\\
&\implies\frac{\alpha_{0}}{2C_{k}'}e^{s-f(s)}=\phi(\tau)\big[\big(\frac{\alpha_{0}}{2C_{k}}(k+1)!-\frac{2\alpha_{1}k(k+1)}{(2\pi)^{2}}+\tilde{C_{k}}\big)(1+\tau)^{k}+\frac{2\alpha_{1}k(1+\tau)^{k-2}}{(2\pi)^{2}}\big].
	\end{split}
\end{equation}
We set 
\begin{equation}
	\label{ D k tau}
	D_{k}(\tau)=	\left(\frac{\alpha_{0}}{2C_{k}}(k+1)!-\frac{2\alpha_{1}k(k+1)}{(2\pi)^{2}}+\tilde{C_{k}}\right)(1+\tau)^{k}+\frac{2\alpha_{1}k(1+\tau)^{k-2}}{(2\pi)^{2}}.
\end{equation}
Assume that $\alpha_{0}>0$ and $D_{k}(\tau)>0$ for $\tau\in[0,2]$.  
 Taking $\log$ on both sides of (\ref{K even equation 1}) and differentiating with respect to $\tau$, we get 
\begin{equation}
	\begin{split}
			&\frac{\big[\phi(\tau)D_{k}(\tau)\big]'}{\phi(\tau)D_{k}(\tau)}=\frac{1-\tau}{\phi(\tau)}\\
		&\phi(\tau)=\frac{\int_{0}^{\tau}(1-t) D_{k}(t)dt}{D_{k}(\tau)}.
	\end{split}
\end{equation}
In the first line we have used (\ref{Legendre transform relation}). We already have that $\phi(0)=0$ and $\phi'(0)=1$. Now if we can show that $\phi(2)=0$, then it also implies $\phi'(2)=-1$. Now $\phi(2)=0$ implies 
\begin{equation}
	\label{phi 2 zero condition in even case}
	\begin{split}
		&\int_{0}^{2}(1-t)D_{k}(t)dt=0\\
		&\implies \frac{\alpha_{0}}{2C_{k}}(k+1)!-\frac{2\alpha_{1}k(k+1)}{(2\pi)^{2}}+\tilde{C_{k}}=-\frac{2\alpha_{1}k}{(2\pi)^{2}}R_{k},
	\end{split}
\end{equation}
where $R_{k}$ is defined by (\ref{ R k}). Now putting (\ref{phi 2 zero condition in even case}) in (\ref{ D k tau}), we get 
\begin{equation}
	\label{D k tau final}
	D_{k}(\tau)=\frac{2\alpha_{1}k(1+\tau)^{k-2}}{(2\pi)^{2}}\left(1-R_{k}(1+\tau)^{2}\right).
\end{equation}
\subsubsection{ k=2 case}
We have already seen that $R_{2}=0$ (\ref{R k calculations}). Hence for $\alpha_{1}>0$, we see that  $D_{2}(\tau)>0$ for all $\tau\in [0,2]$. In this subcase, the condition (\ref{phi 2 zero condition in even case}) becomes  
\begin{equation}
	\frac{3\alpha_{0}}{C_{2}}-\frac{12\alpha_{1}}{(2\pi)^{2}}+\tilde{C_{2}}=0.
\end{equation}
 We can calculate $C_{2}$ and $\tilde{C_{2}}$ explicitly.
Now
\begin{equation}
		C_{2}=52(2\pi)^{3}
\end{equation}
which follows from (\ref{C_{k} expression})
and 
\begin{equation}
	\begin{split}
		&\tilde{C_{2}}=2\alpha_{1}\frac{C_{2}-2\int_{X_{2}}\omega_{\Sigma_{1}}\wedge\omega_{\Sigma_{2}}\wedge d\tau\wedge d\theta}{(2\pi)^{2}\frac{C_{2}}{6}}\\
		&=\frac{12\alpha_{1}}{(2\pi)^{2}}(1-\frac{1}{13}).
	\end{split}
\end{equation}
 We have the following theorem.
\begin{theorem}
	Suppose $\phi:[0,2]\mapsto \mathbb{R}$ is the function $\phi(x)=x-\frac{x^{2}}{2}$. Then by the above construction, there is a solution of the coupled equations satisfying the following conditions-
	\begin{equation}
		\begin{split}
		& 1) \alpha_{0},\alpha_{1}>0\\
		& 2) \alpha_{0}=32\pi\alpha_{1}
		\end{split}
	\end{equation}	  
\end{theorem}
\subsubsection{ k=4 case}
We have already seen that $R_{4}=\frac{5}{46}$ (\ref{R k calculations}). Hence for $\alpha_{1}>0$, we see that  $D_{4}(\tau)>0$ for all $\tau\in [0,2]$.  In this subcase, the condition (\ref{phi 2 zero condition in even case}) becomes  
\begin{equation}
	\frac{60\alpha_{0}}{C_{4}}-\frac{40\alpha_{1}}{(2\pi)^{2}}+\tilde{C_{4}}=-\frac{8\alpha_{1}}{(2\pi)^{2}}\frac{5}{46}.
\end{equation}
 Calculating $C_{4}$ and $\tilde{C}_{4}$ explicitly in this case (using \ref{C_{k} expression}), we get 
\[C_{4}=4!(2\pi)^{5}(3^{5}-1);\ \ \ \tilde{C}_{4}=\frac{40\alpha_{1}}{(2\pi)^{2}}(1-\frac{13}{363}).\]
\begin{theorem}
	Suppose $\phi:[0,2]\rightarrow \mathbb{R}$ is the function $\phi(x)=\frac{(-\frac{40\alpha_{1}}{46} )\int_{0}^{x}(1-t)(1+t)^{4}dt+8\alpha_{1}\int_{0}^{x}(1-t)(1+t)^{2}dt}{(-\frac{40\alpha_{1}}{46} )(1+x)^{4}+8\alpha_{1}(1+x)^{2}}$. Then by the above construction, there is a solution of the coupled equations satisfying the following conditions-
	\begin{equation}
		\begin{split}
			& 1) \alpha_{0},\alpha_{1}>0\\
			& 2) \frac{2\times 235}{3\times 23}\alpha_{1}=\frac{3\alpha_{0}}{4!(2\pi)^{3}}.
		\end{split}
	\end{equation}
\end{theorem}
\subsubsection{ k$\geq$6 case}
\label{k greater than or equals to 6 case}
Using lemma (\ref{Lemma}), we see that (\ref{D k tau final}) can not have a definite sign in $[0,2]$ for any choice of $\alpha_{1}$. Hence in this case, our method does not produce any solutions.
\subsection{Case 2: k odd}
In this case, we take $\gamma_{k}=(k-1)\omega_{\Sigma_{1}}-\sum_{i=2}^{k}\omega_{\Sigma_{i}}$ as our traceless, closed, $(1,1)$ form. One can check that it is indeed traceless and closed. The cohomology class $\frac{\omega+\gamma_{k}}{2\pi}$ is integral. Hence there exists a Hermitian line bundle $(L_{k},H_{k})$ with curvature $\sqrt{-1}F_{H_{k}}=\frac{\omega+\gamma_{k}}{2\pi}$.\par
Now $k!\sqrt{-1}\bigwedge_{\omega}F_{H_{k}}\frac{\omega^{k+1}}{(k+1)!}=\sqrt{-1}F_{H_{k}}\wedge\omega^{k}=\frac{(\omega+\gamma_{k})}{2\pi}\wedge\omega^{k}=\frac{1}{2\pi}\omega^{k+1}$. So $\sqrt{-1}\bigwedge_{\omega}F_{H_{k}}=\frac{k+1}{2\pi}$. Hence, the first equation is met with $\lambda=\frac{k+1}{2\pi}$.\par
Now 
\begin{equation}
	\label{expression for alpha square in k odd case}
	\gamma_{k}^{2}=-2(k-1)\sum_{i=2}^{k}\omega_{\Sigma_{1}}\wedge\omega_{\Sigma_{i}}+2\sum_{2\leq i<j\leq k}\omega_{\Sigma_{i}}\wedge\omega_{\Sigma_{j}}
\end{equation}
Using (\ref{expression for alpha square in k odd case}) and (\ref{expression for omega k-1}) , we get that 
\begin{equation}
	\label{expression for alpha square omega to the power k-1}
\begin{split}
	\gamma_{k}^{2}\wedge\omega^{k-1}&=\left(-2(k-1)\sum_{i=2}^{k}\omega_{\Sigma_{1}}\wedge\omega_{\Sigma_{i}}+2\sum_{2\leq i<j\leq k}\omega_{\Sigma_{i}}\wedge\omega_{\Sigma_{j}}\right)\wedge\\
	&(k-1)!(1+\tau)^{k-2}\phi(\tau)\sqrt{-1}\frac{dw\wedge d\bar{w}}{\lvert w\rvert^{2}}\wedge\sum_{l<m}\omega_{\Sigma_{1}}\wedge\dots\wedge\widehat{\omega_{\Sigma_{l}}}\wedge\dots\wedge\widehat{\omega_{\Sigma_{m}}}\wedge\dots\wedge\omega_{\Sigma_{k}}\\
	&=-k(k-1)(k-1)!(1+\tau)^{k-2}\phi(\tau)\omega_{\Sigma_{1}}\wedge\dots\wedge\omega_{\Sigma_{k}}\wedge\sqrt{-1}\frac{dw\wedge d\bar{w}}{\lvert w\rvert^{2}}
\end{split}	
\end{equation}
Putting (\ref{expression for omega to the power k+1}) and (\ref{expression for alpha square omega to the power k-1}) in (\ref{second Chern character }) , we get 
\begin{equation}
	\label{second Chern character in odd case}
\frac{2}{(k-1)!}F_{H}\wedge F_{H}\wedge \omega^{k-1}=-\frac{2\phi(\tau)}{(2\pi)^{2}(k-1)!}[(k+1)!(1+\tau)^{k}-k(k-1)(k-1)!(1+\tau)^{k-2}]\omega_{\Sigma_{1}}\wedge\dots\wedge\omega_{\Sigma_{k}}\wedge\frac{\sqrt{-1}dw\wedge d\bar{w}}{\lvert w\rvert^{2}}\\
\end{equation}
Now putting (\ref{expression for Omega h}) , (\ref{expression for omega to the power k+1}) , (\ref{second Chern character in odd case}), (\ref{expression of volume form}) and $\int_{X}\Omega_{h}=C_{k}'$, $\int_{X}\omega^{k+1}=C_{k}$   in the second equation of (\ref{Coupled equations}) , we get 
\begin{equation}
	\begin{split}
		\label{equation in odd case}
		&\frac{\alpha_{0}}{2}[\frac{e^{s-f}}{C_{k}'}-\frac{(k+1)!(1+\tau)^{k}\phi(\tau)}{C_{k}}]+\frac{2\alpha_{1}\phi(\tau)}{(2\pi)^{2}(k-1)!}[(k+1)!(1+\tau)^{k}-k(k-1)(k-1)!(1+\tau)^{k-2}]\\
		&=\tilde{C_{k}}(1+\tau)^{k}\phi(\tau)\\
		&\implies \phi(\tau)\big[\big(\frac{\alpha_{0}(k+1)!}{2C_{k}}-\frac{2\alpha_{1}k(k+1)}{(2\pi)^{2}}+\tilde{C_{k}}\big)(1+\tau)^{k}+\frac{2\alpha_{1}k(k-1)}{(2\pi)^{2}}(1+\tau)^{k-2}\big]=\frac{\alpha_{0}}{2C_{k}'}e^{s-f(s)}
	\end{split}
\end{equation}
 We set 
\begin{equation}
\label{positivity in odd case}
P_{k}(\tau)=\big(\frac{\alpha_{0}(k+1)!}{2C_{k}}-\frac{2\alpha_{1}k(k+1)}{(2\pi)^{2}}+\tilde{C_{k}}\big)(1+\tau)^{k}+\frac{2\alpha_{1}k(k-1)}{(2\pi)^{2}}(1+\tau)^{k-2}.
\end{equation}
Assume that $\alpha_{0}>0$ and $P_{k}(\tau)>0$ for $\tau \in [0,2]$. Then taking $\log$ on both side of (\ref{equation in odd case}) and differentiating with respect to $\tau$, we get 
\begin{equation}
	\begin{split}
		&\frac{\big[\phi(\tau)P_{k}(\tau)\big]'}{\phi(\tau)P_{k}(\tau)}=\frac{1-\tau}{\phi(\tau)}\\
		&\phi(\tau)=\frac{\int_{0}^{\tau}P_{k}(t)(1-t)dt}{P_{k}(\tau)},
	\end{split}
\end{equation}
where we have used (\ref{Legendre transform relation}). We already have that $\phi(0)=0$ and $\phi'(0)=1$. Now if we can show that $\phi(2)=0$, then it also implies $\phi'(2)=-1$. Now $\phi(2)=0$ implies 
 \begin{equation}
 	\label{phi 2 zero odd case}
 	\begin{split}
 		&\int_{0}^{2}(1-t)P_{k}(t)dt=0\\
 		&\implies \frac{\alpha_{0}(k+1)!}{2C_{k}}-\frac{2\alpha_{1}k(k+1)}{(2\pi)^{2}}+\tilde{C_{k}}=-\frac{2\alpha_{1}k(k-1)}{(2\pi)^{2}}R_{k},
 	\end{split}
 \end{equation}
where $R_{k}$ is given by (\ref{ R k}). Now putting (\ref{phi 2 zero odd case}) in (\ref{positivity in odd case}), we have 
\begin{equation}
	\label{P k tau final}
	P_{k}(\tau)=\frac{2\alpha_{1}k(k-1)(1+\tau)^{k-2}}{(2\pi)^{2}}\left(1-R_{k}(1+\tau)^{2}\right)
\end{equation}
\subsubsection{ k=3 case}
We have already seen that $R_{3}=\frac{5}{63}$. Hence for $\alpha_{1}>0$, we see that $P_{k}(\tau)>0$ for all $\tau\in[0,2]$. In this subcase, the condition (\ref{phi 2 zero odd case}) becomes  
\begin{equation}
	\frac{12\alpha_{0}}{C_{3}}-\frac{24\alpha_{1}}{(2\pi)^{2}}+\tilde{C_{3}}=-\frac{12\alpha_{1}}{(2\pi)^{2}}\times\frac{5}{63}.
\end{equation}
 In this case, we can calculate $C_{3}$ (using \ref{C_{k} expression}) and $\tilde{C_{3}}$ explicitly. Indeed,
\[C_{3}=3!(2\pi)^{4}(3^4-1);\ \ \ \tilde{C_{3}}=\frac{24\alpha_{1}}{(2\pi)^{2}}(1-\frac{1}{10}).\]
 So we have the following theorem.
\begin{theorem}
	Suppose $\phi:[0,2]\rightarrow \mathbb{R}$ is the function $\phi(x)=\frac{(-\frac{60\alpha_{1}}{63} )\int_{0}^{x}(1-t)(1+t)^{3}dt+12\alpha_{1}\int_{0}^{x}(1-t)(1+t)dt}{(-\frac{60\alpha_{1}}{63} )(1+x)^{3}+12\alpha_{1}(1+x)}$. Then by the above construction, there is a solution of the coupled equations satisfying the following conditions-
	\begin{equation}
		\begin{split}
			& 1) \alpha_{0},\alpha_{1}>0\\
			& 2) \frac{\alpha_{0}}{32(2\pi)^{2}}=\frac{38}{21}\alpha_{1}.
		\end{split}
	\end{equation}
\end{theorem}
\subsubsection{ $k\geq$5 case}
\label{k greater than or equals to 5 case}
Using lemma (\ref{Lemma}), we see that (\ref{P k tau final}) can not have a definite sign in $[0,2]$ for any choice of $\alpha_{1}$. Hence in this case, our method does not produce any solutions.
\begin{remark}
	 One can take the traceless closed $(1,1)$ form as $\gamma_{k}=\sum_{i=1}^{i=k}c_{i}\omega_{\Sigma_{i}}$, where $\sum_{i=1}^{i=k}c_{i}=0$ and $c_{i}$ are integers.
\end{remark}

Department of Mathematics, Indian Institute of Science, Bangalore, India - $560012$\\
E-mail address: \textit{kartickghosh@iisc.ac.in}
\end{document}